%% file: mpleGenHAL2.tex
\documentclass[preprint]{article}

\usepackage{graphicx}
\usepackage[utf8]{inputenc}
\usepackage[T1]{fontenc}
\usepackage{aeguill}
\usepackage{amssymb}
\usepackage[utf8]{inputenc}
\usepackage[T1]{fontenc}
\usepackage{aeguill}
\usepackage{amssymb}
\usepackage{graphicx}
\usepackage[a4paper,pdftex=true]{geometry}

\usepackage{natbib}
\usepackage{bbm}
\usepackage{amsmath}
\usepackage{hyperref}
\usepackage[a4paper,pdftex=true]{geometry}
\usepackage{fancyvrb}

\newtheorem{theorem}{Theorem}
\newtheorem{lemma}[theorem]{Lemma}
\newtheorem{corollary}[theorem]{Corollary}

\newtheorem{proposition}[theorem]{Proposition}

\newtheorem{definition}{Definition}

\input cqlsInclude

\input testInclude

\renewcommand{\RR}[1][2]{\mathbb{R}^{#1}}
\renewcommand{\ZZ}[1][2]{\mathbb{Z}^{#1}}
\renewcommand{\NN}[1][]{\mathbb{N}^{#1}}

\newcommand{\cSp}[1][]{{\Omega}#1
}

\newcommand{\cSptt}{\overline{\Omega}
}

\newcommand{\ism}[2][\Lambda]{
\int_{#1} #2 dx
}

\newcommand{\cD}[1][]{
\left\lceil{#1}\right\rceil
}

\newcommand{\VE}{V}

\newcommand{\Par}{\theta}                    
\newcommand{\ParV}{\Vect{\theta}}            

\newcommand{\ParVT}{\Vect{\theta}^\star}     
\newcommand{\SpPar}{\Vect{\Theta}}           
\newcommand{\Vois}{\mathcal{{V}}}

\newcommand{\PL}[3]{PL_{#1} \left( #2;#3 \right)}
\newcommand{\LPL}[3]{LPL_{#1} \left( #2;#3 \right)}
\newcommand{\Un}[2][\varphi]{\UnE\left( #1;#2 \right)}
\newcommand{\UnE}{U_n
}
\newcommand{\ddU}[1][\ParV]{\Mat{U}^{(2)}(#1)
}

\newcommand{\Dom}[1][0]{\Delta_{#1}}

\newcommand{\VI}[2]{\VE\left(#1|#2\right)}

\newcommand{\VIPar}[3]{\VE \left(#1|#2  ; #3 \right)}

\newcommand{\dVIPar}[4]{\frac{\partial \VE}{\partial \Par_{#1}} \left(#2|#3  ; #4 \right)}
\newcommand{\ddVIPar}[5]{\frac{\partial^2 \VE}{\partial \Par_{#1} \partial \Par_{#2}  } \left(#3|#4  ; #5 \right) }

\newcommand{\Estn}[1]{ \widehat{#1}_n }

\newcommand{\Dt}{\widetilde{D}}

\newcommand{\Ltt}{\overline{\Lambda}}

\newcommand{\BDt}[2][i]{\ensuremath{\mathbbm{B}\left(#1,#2\right)}}

\begin{document}

\title{Asymptotic properties \\ of the maximum pseudo-likelihood estimator \\ for stationary Gibbs point processes  including \\ the Lennard-Jones model}

  \author{Jean-François Coeurjolly$^{1,2}$ and Rémy Drouilhet$^{1}$ \\ 
$^1$ LJK, Grenoble University, France \\
$^2$ GIPSA-lab, Grenoble University, France}



\maketitle

\begin{center}\begin{minipage}{13cm}{\centerline{\bf Abstract} 
This paper presents asymptotic properties of the maximum pseudo-likelihood estimator of a vector $\Vect{\theta}$ parameterizing a stationary Gibbs point process. Sufficient conditions, expressed in terms of the local energy function defining a Gibbs point process, to establish strong consistency and asymptotic normality results of this estimator depending on a single realization, are presented. 
These results are general enough to no longer require the local stability and the linearity in terms of the parameters of the local energy function. 
We consider characteristic examples of such models, the Lennard-Jones and the finite range Lennard-Jones models. We show that the different assumptions ensuring the consistency are satisfied for both models whereas the assumptions ensuring the asymptotic normality are fulfilled only for the finite range Lennard-Jones model.\\
}\end{minipage}\end{center}

\noindent {\bf Keywords}: Stationary Gibbs point processes, maximum pseudo-likelihood estimator, Lennard-Jones model.

\section{Introduction}

These last years, much attention has been paid to spatial point pattern data, and especially to models and methodologies for fitting them, see \cite{R-Mol08} for a recent overview of this topic and  \cite{B-DalVer88}, \cite{B-StoKenMecRus87} \cite{B-MolWaa03} or \cite{B-IllPenSto08} for  more general information. For spatial point pattern data, the reference model is the Poisson point process modelling a random configuration of points with no interaction between points. In particular, this leads to the independence of any two random sub-configurations lying in two non-overlapping domains. 
A way to introduce dependence is to consider the class of Gibbs models. In a bounded domain, a Gibbs point process is defined through its probability measure having a Radon-Nykodym derivative with respect to a Poisson point process measure proportional to $e^{-V(\varphi)}$ where $V(\varphi)$ corresponds to the energy function (i.e. a cost function expressed in terms of interactions) of the configuration of points $\varphi$. The definition of Gibbs models in $\RR[d]$ is essential  when dealing with asymptotic properties of estimators based on a point process observed in a domain aimed at converging towards $\RR[d]$. The extension of this definition is not so straightforward. The probability measure of a Gibbs point process in $\RR[d]$ has to be defined by specifying its conditional density (indirectly expressed in terms of the energy function $V(\varphi)$), see  {\it e.g.} \cite{B-Pre76} or Section~\ref{sec-background} for more details. 

The class of Gibbs point processes is extremely rich. The energy function can penalize points, pairs or triplets of points (see {\it e.g.} \cite{A-BadTur00}). More sophisticated models can also be obtained by considering interactions based on the Delaunay or the $k-$nearest neighbor graphs (\cite{BBD1,BBD2}), Voronoï tessellations (\cite{A-DerLav09}) or random sets (\cite{A-KenLieBad99}, \cite{A-Der09}).

Following the definition of a parametric Gibbs point process, the natural question of efficiently estimating the parameters arises. Many proposals have tried to estimate the energy function from an available point pattern data. The most well-known method is the use of the likelihood function, see {\it e.g. } \cite{B-MolWaa03} and the references therein. The main drawback of this approach is that the likelihood function contains an unknown scaling factor whose value depends on the parameters. This parametric normalizing constant is difficult to calculate from a practical point of view. From a theoretical one, it also makes asymptotic results more complicated to obtain. An alternative approach relies on the use of the pseudo-likelihood function. The idea originated from \cite{Besag74} in the study of lattice processes. \cite{Besag82} further considered this method for pairwise interaction point processes, and \cite{Jensen91} extended the definition of the pseudo-likelihood function to the general class of marked Gibbs point processes. The construction of the pseudo-likelihood function is based on the conditional densities which spare the computation of the scaling factor.

Our paper deals with asymptotic properties of the maximum pseudo-likelihood estimator. In order to underline our theoretical improvements, let us discuss the two main different papers discussing this topic:
\begin{itemize}
\item In \cite{A-BilCoeDro08}, we obtain consistency and asymptotic normality for exponential family models of Gibbs point processes, that is, on models with energy functions that are linear in terms of the parameters. Moreover, we concentrate on models such that the local energy function is local and stable. The locality of the local energy expresses that the energy to insert a point $x$ into $\varphi$, that is, $V(x|\varphi)=V(\varphi\cup x)-V(\varphi)$, depends only on the points of $\varphi$ falling into some ball with a fixed radius whereas the stability of the local energy (property referred as the local stability) asserts that $V(x|\varphi)$ is bounded from below by a finite negative constant. The paper \cite{A-BilCoeDro08} extends several papers (\cite{Jensen91}, \cite{Jensen94}) and includes a large class of examples of practical interest: area-interaction point process, Multi-Strauss marked point process based on the complete graph or the $k$-nearest-neighbors graph, or the Geyer's triplet point process to name a few. 
\item Another work has been undertaken by Mase. The consistency for non necessarily stable local energy functions (actually
for superstable and lower regular ones introduced by \cite{Ruelle70}) is obtained in \cite{Mase95} for specific models with only two parameters -the chemical potential and the inverse temperature-  which can be viewed as particular exponential family models. \cite{Mase00} extended his work to the context of marked point processes and provided asymptotic normality by adding the assumption of finite range. 
\end{itemize}

Based on this literature, the main goal of this paper is to derive asymptotic properties similar to the ones presented before (consistency and asymptotic normality) but in a more general framework. We provide asymptotic results for general Gibbs point processes with {\it non (necessarily) linear} and {\it non (necessarily) stable} local energy functions. The characteristic example we have in mind is the {\it Lennard-Jones model}. This model, from statistical physics, is a stationary pairwise interaction Gibbs point process where the local energy to insert a point $x$ into a configuration $\varphi$ is parameterized as follows: for $\ParV=(\theta_1,\theta_2,\theta_3)\in \RR[3]$ with $\theta_2,\theta_3>0$ 
$$
V^{LJ}\left( x|\varphi; \ParV \right) 
:= \theta_1 +  4 \theta_2 \sum_{y\in \varphi}\left( \left(\frac{\theta_3}{\| y-x\|}\right)^{12}
-\left(\frac{\theta_3}{\|y-x\|}\right)^6 \right).
$$
Let us notice that \cite{Mase95} could only propose the estimation of $\theta_1$ and $\theta_2$ with known $\theta_3$.
The Lennard-Jones model is of great interest from several points of view. From a physical point of view, this model arises when theoretically modelling a pair of neutral atoms or molecules subject to two distinct forces in the limit of large separation and small separation: an attractive force at long ranges (van der Waals force, or dispersion force) and a repulsive force at short ranges (the result of overlapping electron orbitals, referred to as a Pauli repulsion from the Pauli exclusion principle). In this literature, the parameters $\theta_2$ and $\theta_3$ are often referred to as the depth potential and the (finite) distance at which the interparticle potential is zero. From a probabilistic point of view, this model constitutes the main example of superstable, regular and lower regular energies studied in \cite{Ruelle70} where the author proves the existence of ergodic measures for such models. Finally, from a statistical point of view, this model has been considered by several authors, see {\it e.g.} \cite{A-OgaTan81}, \cite{A-GouSarGra96} for fitting spatial point patterns arising in forestry. In particular, let us note that, in \cite{A-GouSarGra96}, the model is fitted by using the maximum pseudo-likelihood method. As the authors do not endeavour to justify the theoretical performances of the procedure, the result proposed in Section~\ref{sec-LJ} of this paper fills this gap.

The rest of the paper is organized as follows. Section~\ref{sec-background} introduces some background and notation on Gibbs point processes (general definitions, examples). The maximum pseudo-likelihood method and asymptotic results of the derived estimator are proposed in Section~\ref{sec-results}. For general Gibbs point processes, sufficient conditions, expressed in terms of the local energy function to establish strong consistency and asymptotic normality results of this estimator are presented. While no general condition on the model is assumed to obtain the consistency, the characteristic finite range of the local energy function is required to establish the asymptotic normality.
For the sake of simplicity, Section~\ref{sec-results} (and the resulting proofs) would concentrate on non-marked Gibbs point processes. However, as we have shown in our seminal paper \cite{A-BilCoeDro08}, no real mathematical difficulty occurs with the introduction of marks. Section~\ref{sec-LJ} focuses on the Lennard-Jones model. We show that the general assumptions described in Section~\ref{sec-results} are fulfilled for this model. Proofs have been postponed until Section~\ref{sec-proofs}.


\section{Background and notation} \label{sec-background}

For the sake of simplicity, we consider Gibbs point processes in dimension $d=2$. 

\subsection{General notation, configuration space}
Subregions of $\RR[2]$ will typically be denoted by $\Lambda$ or $\Delta$ and will always be assumed to be Borel with positive Lebesgue measure. We write $\Lambda \Subset \RR[2]$ if $\Lambda$ is bounded. $\Lambda^c$ denotes the complementary set of $\Lambda$ inside $\RR[2]$.
The notation $|.|$ will be used without ambiguity for different kind of objects. For a countable set $\mathcal J$, $|\mathcal J|$ represents the number of elements belonging to $\mathcal J$; For $\Lambda\Subset\RR[2]$, $|\Lambda|$ is  the volume of $\Lambda$; For a vector $x \in \RR[2]$, $|x|$ corresponds to its uniform norm while $\| x\|$ is simply its euclidean norm.
For all $ x\in\RR[2],\rho>0$ and $i\in\ZZ[2]$, let $\mathcal{B}( x,\rho):=\{ y \in \RR[2],\ | y- x|<\rho\}$ and $\mathbbm{B}(i,\rho) := \mathcal{B}(i,\rho) \cap \ZZ[2]$.

A configuration is a subset $\varphi$ of $\RR[2]$ which is locally finite in that $\varphi_\Lambda:=\varphi \cap \Lambda$ has finite cardinality $N_\Lambda(\varphi):=|\varphi_\Lambda|$ for all $\Lambda \Subset \RR[2]$. The space $\Omega$ of all configurations is equipped with the $\sigma$-algebra $\mathcal{F}$ that is generated by the counting variables $N_\Lambda(\varphi)$ with $\Lambda \Subset \RR[2]$. Finally, let $T=\left( \tau_x \right)_{x\in \RR[2]}$ be the shift group, where $\tau_x:\Omega\to \Omega$ is the translation by the vector $-x \in \RR[2]$. 

\subsection{Gibbs point processes}

Our results will be expressed for general stationary Gibbs point processes. Since we are interested in asymptotic properties, we have to consider these point processes acting on the infinite volume $\RR[2]$. Let us briefly recall their definition.

A point process $\Phi$ is a $\Omega$-valued random variable, with probability distribution $P$ on $(\Omega,\mathcal F)$. 
The most prominent point process is the (homogeneous) Poisson process with intensity $z>0$. Recall that its probability measure $\pi^z$ is the unique probability measure on $(\Omega,\mathcal{F})$ such that the following holds for $\Lambda \Subset \RR[2]$: $(i)$ $N_\Lambda$ is Poisson distributed with parameter $z|\Lambda|$, and $(ii)$ conditionally to $N_\Lambda=n$, the $n$ points in $\Lambda$ are independent with uniform distribution on $\Lambda$, for each interger $n\geq 1$. For $\Lambda\Subset \RR[2]$, let us denote  by $\pi^z_\Lambda$ the marginal probability measure in $\Lambda$ of the Poisson process with intensity $z$.

Let $\ParV \in \RR[p]$ (for some $p\geq 1$). For any $\Lambda\Subset \RR[2]$, let us consider the parametric function $V_{\Lambda}(.;\ParV)$ from $\Omega$ into $\RR[]\cup\{+\infty\}$. From a physical point of view, $V_{\Lambda}(\varphi;\ParV)$ is the energy of $\varphi_{\Lambda}$ in $\Lambda$ given the outside configuration $\varphi_{\Lambda^c}$. 

In this article, we focus on stationary point processes on $\RR[2]$, i.e. with $T$-invariant probability measure. For any $\Lambda\Subset \RR[2]$, we therefore consider $V_{\Lambda}(.;\ParV))$ to be $T$-invariant, i.e. $V_{\Lambda}(\tau_x \varphi;\ParV)=V_{\Lambda}(\varphi;\ParV)$ for any $x\in\RR[2]$. Furthermore, we assume that the family of energies is hereditary, which means that for any $\Lambda\Subset \RR[2]$, $\varphi\in\Omega$, and $x \in\Lambda$:
$V_{\Lambda}(\varphi;\ParV))=+\infty \Rightarrow V_{\Lambda}(\varphi \cup\{x\};\ParV))=+\infty.$

In such a context, a Gibbs measure is usually defined as follows (see \cite{B-Pre76}).

 \begin{definition} \label{def-Gibbs} A probability measure $P_{\ParV}$ on $\Omega$  is a Gibbs measure for the  family of energies $(V_{\Lambda}(.;\ParV))_{\Lambda\Subset \RR[2]}$ if for every $\Lambda\Subset \RR[2]$, for $P_{\ParV}$-almost every outside configuration $\varphi_{\Lambda^c}$, the law of $P_{\ParV}$ given $\varphi_{\Lambda^c}$ admits the following density with respect to $\pi_\Lambda^{z}$:
 $$f_{\Lambda}(\varphi_{\Lambda}|\varphi_{\Lambda^c};\ParV)=\frac{1}{Z_{\Lambda}(\varphi_{\Lambda^c};\ParV)}e^{-V_{\Lambda}(\varphi;\ParV)},$$
where $Z_{\Lambda}(\varphi_{\Lambda^c};\ParV):=\int_{\Omega_\Lambda} e^{-V_{\Lambda}(\varphi_\Lambda\cup \varphi_{\Lambda^c};\ParV)} \pi_\Lambda^z(d\varphi_\Lambda)$ is called the partition function.
\end{definition}
Without loss of generality, the intensity of the Poisson process, $z$ is fixed to~1 and we simply write $\pi$ and $\pi_\Lambda$ in place of $\pi^1$ and $\pi_\Lambda^1$.
In the previous definition, we implicitly assume the consistency of the family $(f_{\Lambda}(.|.;\ParV))_{\Lambda\Subset \RR[2]}$: for any $\Delta\subset\Lambda\Subset\RR[2]$
$$
f_{\Delta}(\varphi_{\Delta}|\varphi_{\Delta^c};\ParV)= \frac{f_{\Lambda}(\varphi_{\Delta}\cup \varphi_{\Lambda\setminus \Delta}|\varphi_{\Lambda^c};\ParV)}{f_{\Lambda}(\varphi_{\Lambda\setminus \Delta}|\varphi_{\Lambda^c};\ParV)}=\frac{f_{\Lambda}(\varphi_{\Delta}\cup \varphi_{\Lambda\setminus \Delta}|\varphi_{\Lambda^c};\ParV)}{\int_{\Omega_\Delta}f_{\Lambda}(\psi_\Delta\cup \varphi_{\Lambda\setminus \Delta}|\varphi_{\Lambda^c};\ParV) \pi_\Delta(d\psi_\Delta)}.
$$
A sufficient condition to directly fulfill this basic ingredient is to assume the compatibility of the family  $(V_\Lambda(.))_{\Lambda\Subset \RR[2]}$: for every $\Delta\subset\Lambda \Subset\RR[2]$, the function  $\varphi \to V_{\Lambda}(\varphi;\ParV)-V_{\Delta}(\varphi;\ParV)$ from $\Omega$ into $\RR[]\cup\{+\infty\}$ is measurable and only depends on $\varphi_{\Lambda^c}$.

The existence of a Gibbs measure on $\Omega$ which satisfies these conditional specifications is a difficult issue. We refer the interested reader to \cite{B-Rue69,B-Pre76,A-BerBilDro99,A-Der05,A-DerDroGeo09} for the technical and mathematical development of the existence problem. The minimal assumption of our paper is then:
\begin{list}{}{}
\item \textbf{[Mod-E]}: Our data consist in the realization of a point process $\Phi$ with Gibbs measure $P_{\ParVT}$, where $\ParVT \in \mathring{\SpPar}$, $\SpPar$ is a compact subset of $\RR[p]$ and, for any $\ParV \in\SpPar$, there exists a stationary Gibbs measure $P_{\ParV}$ for the  family $(V_{\Lambda}(.;\ParV))_{\Lambda\Subset \RR[2]}$.
\end{list}

In the rest of this paper, the reader has mainly to keep in mind the concept of  local energy defined as the energy required to insert a point $x$ into the configuration $\varphi$ and expressed for any $\Lambda\ni x$ by 
$$V^{}\left( x|\varphi; \ParV \right):=V_{\Lambda}(\varphi \cup\{x\})-V_{\Lambda}(\varphi).$$
From the compatibility of the family of energies, the local energy does not depend on $\Lambda$.

Our asymptotic normality result will require the following locality property assumption. 


\begin{list}{}{}
\item \textbf{[Mod-L]}: There exists $D\geq 0$ such that for all $\varphi \in  \Omega$
$$\VIPar{0}{\varphi}{\ParV} = \VIPar{0}{\varphi_{\mathcal B(0,D)}}{\ParV}.$$
\end{list}

\subsection{Example : Lennard-Jones models} \label{sec-LJtype}

Let us present the main example studied in this paper. We call LJ-type model the stationary pairwise interaction point process defined for some $D\in ]0,+\infty]$ by
\[
V^{LJ}_\Lambda\left(\varphi; \ParV \right) :=\theta_1|\varphi_\Lambda| +  H^{LJ}_\Lambda\left(\varphi; \ParV \right) \mbox{ with }
 H^{LJ}_\Lambda\left(\varphi; \ParV \right):=\sum_{\begin{subarray}{c}x_1\in \varphi_\Lambda\\x_2\in\varphi_{\Lambda^c}\end{subarray}} g^{LJ}(||x_1-x_2||;\ParV)
\]
and \[
g^{LJ}(r;\ParV):= 4 \theta_2\left( \left(\frac{\theta_3}r\right)^{12}-\left(\frac{\theta_3}r\right)^{6}\right)\mathbf{1}_{[0,D]}(r).
\]
As a direct consequence, the local energy function is expressed as
\[
 V^{LJ}\left( x|\varphi; \ParV \right) :=\theta_1 +  H^{LJ}\left( x|\varphi; \ParV \right) \mbox{ with }
 H^{LJ}\left( x|\varphi; \ParV \right):=\sum_{y \in \varphi} g^{LJ}(||x-y||;\ParV).
\]
where $\ParV=(\theta_1,\theta_2,\theta_3)\in \RR[]\times(\RR[+])^2$. The cases  $D=+\infty$ and $D<+\infty$ respectively correpond to the Lennard-Jones model (briefly presented in the introduction) and the Lennard-Jones model with finite range.

\cite{Ruelle70} has proved the existence of an ergodic measure for superstable, regular and lower regular potentials. The Lennard-Jones model (including the finite range one) is known to be the characteristic example of such a family of models for which Ruelle managed to prove the existence of ergodic measures for any $\ParV\in \RR[]\times(\RR[+])^2$.  
In order to ensure \textbf{[Mod-E]}, it is required to assume that $\theta_2^\star,\theta_3^\star>0$. 
Finally, \textbf{[Mod-L]} is satisfied for the LJ-type model with $D<+\infty$ since the parameter $D$ corresponds for pairwise interaction point processes to the range of the Gibbs point process.


\section{Asymptotic results of the Maximum pseudo-likelihood estimator} \label{sec-results}

\subsection{Maximum pseudo-likelihood method}

The idea of maximum pseudo-likelihood is due to \cite{Besag74} who first introduced the concept for Markov random fields in order to avoid the normalizing constant. This work was then widely extended and \cite{Jensen91} (Theorem 2.2) obtained a general expression for Gibbs point processes. Using our notation and up to a scalar factor, the pseudo-likelihood defined for a configuration $\varphi$ and a domain of observation $\Lambda$ is denoted by $\PL{\Lambda}{\varphi}{\ParV}$ and given by
\begin{equation} \label{defPLpap}
\PL{\Lambda}{\varphi}{\ParV} = \exp\left( - \ism[\Lambda]{ e^{ - \VIPar{x}{\varphi}{\ParV} }  }\right)
  \prod_{x\in \varphi_\Lambda} e^{ -\VIPar{x}{\varphi \setminus x}{\ParV} }.
\end{equation}
It is more convenient to define and work with the log-pseudo-likelihood, denoted by $\LPL{\Lambda}{\varphi}{\ParV}$
\begin{equation} \label{logPL}
\LPL{\Lambda}{\varphi}{\ParV}=
-  \ism[\Lambda]{ e^{ - \VIPar{x}{\varphi}{\ParV} } }  - \sum_{x\in\varphi_\Lambda}  \VIPar{x}{\varphi \setminus x}{\ParV} .
\end{equation}

The point process is assumed to be observed in a domain $\Lambda_n\oplus\widetilde{D}= \cup_{x \in \Lambda_n}\mathcal{B}(x,\widetilde{D})$ for some $\widetilde{D}<+\infty$. For the asymptotic normality result, it is also assumed that $\widetilde{D}\geq D$ and that $\Lambda_n\subset \RR[2]$ can be decomposed into $\cup_{i \in I_n} \Dom[i]$ 
 where $I_n= \BDt[0]{n}$ and for $i \in \ZZ[2]$,
$\Dom[i]=\Dom[i](\Dt)$ is the square centered at $i$ with side-length $\Dt$. As a consequence, as $n \to +\infty$, $\Lambda_n \to \RR[2]$ such that $|\Lambda_n|\to +\infty$ and ${\displaystyle \frac{|\partial \Lambda_n|}{|\Lambda_n|}\to 0}$.

Define for any configuration $\varphi$, $\Un{\ParV} = - \frac1{|\Lambda_n|} \LPL{\Lambda_n}{\varphi}{\ParV}$. The maximum pseudo-likelihood estimate (MPLE), denoted by $\Estn{\ParV}(\varphi)$, is then defined by
$$
\Estn{\ParV}(\varphi) =  \mathop{\arg\max}_{\ParV \in \SpPar} \; \LPL{\Lambda_n}{\varphi}{\ParV} = \mathop{\arg\min}_{\ParV\in \SpPar} \Un{\ParV}.
$$
The following basic notation are introduced: for $j,k=1,\dots,p$ and $\Lambda \Subset \RR[2]$
\begin{itemize}
\item Gradient vector of $\UnE$: $\ensuremath{\Vect{U}_n^{(1)}(\varphi;\ParV)}:=-|\Lambda_n|^{-1} \ensuremath{\Vect{LPL}_{\Lambda_n}^{(1)} \left( \varphi ; \ParV  \right)}$ where
$$
\left(\ensuremath{\Vect{LPL}_{\Lambda}^{(1)} \left( \varphi ; \ParV  \right)}\right)_j=
\ism[\Lambda]{  {\frac{\partial V^{}}{\partial \theta_{j}}}  \left( x|\varphi; \ParV \right)  e^{ - \VIPar{x}{\varphi}{\ParV} } }  - \sum_{x\in\varphi_\Lambda}     {\frac{\partial V^{}}{\partial \theta_{j}}}  \left( x|\varphi \setminus x; \ParV \right). 
$$
\item Hessian matrix of $\UnE$: $\ensuremath{\Mat{U}_n^{(2)}(\varphi;\ParV)}:=-|\Lambda_n|^{-1} \ensuremath{\Mat{LPL}_{\Lambda_n}^{(2)} \left( \varphi ; \ParV  \right)}$
\begin{eqnarray*}
\left(\ensuremath{\Mat{LPL}_{\Lambda}^{(2)} \left( \varphi ; \ParV  \right)}\right)_{j,k}&=&
\ism{  \left(
\frac{\partial^2 V^{}}{\partial \theta_{j}\partial \theta_{k} }  \left( x|\varphi; \ParV \right)
- \frac{\partial V}{\partial \theta_j} \left( x|\varphi; \ParV \right)\frac{\partial V}{\partial \theta_k} \left( x|\varphi; \ParV \right)
\right)
e^{ - \VIPar{x}{\varphi}{\ParV} } }  \\
&&+ \sum_{x \in \varphi_\Lambda}   \frac{\partial V}{\partial \theta_j} \left( x|\varphi\setminus x; \ParV \right)\frac{\partial V}{\partial \theta_k} \left( x|\varphi \setminus x; \ParV \right).
\end{eqnarray*}

\end{itemize}
Finally, note that from the decomposition of the observation domain $\Lambda_n$, one has
$$
\ensuremath{\Vect{U}_n^{(1)}(\varphi;\ParV)} = -|\Lambda_n|^{-1}  \sum_{i \in I_n} \ensuremath{\Vect{LPL}_{\Delta_i}^{(1)} \left( \varphi ; \ParV  \right)} 
\mbox{ and } 
\ensuremath{\Mat{U}_n^{(2)}(\varphi;\ParV)} = -|\Lambda_n|^{-1}  \sum_{i \in I_n} \ensuremath{\Mat{LPL}_{\Delta_i}^{(2)} \left( \varphi ; \ParV  \right)}. 
$$

\subsection{Consistency of the MPLE}

The assumption \textbf{[C]} gathers the  following four assumptions:
\begin{itemize}
\item[\textbf{[C1]}] For all $\ParV \in \SpPar$,
$$\Esp\left( e^{- \VIPar{0}{\Phi}{\ParV}} \right) <+\infty \quad \mbox{ and } \quad
\Esp\left( \left|\VIPar{0}{\Phi}{\ParV} \right| e^{- \VIPar{0}{\Phi}{\ParVT}}\right) <+\infty.$$
\item[\textbf{[C2]}] 
Identifiability condition~: there exists 
$A_1,\ldots,A_{\ell}$, $\ell\geq p$ events in $\Omega$ such that:
\begin{itemize}
\item the $\ell$ events $A_i$ are disjoint and satisfy $P_{\ParVT}(B_i) >0$ 
\item for all $\left(\varphi_1,\ldots,\varphi_\ell\right) \in A_1 \times\cdots \times A_{\ell}$
$$ \left\{\begin{array}{l}
D(0 |\varphi_i;\ParV)=0\\ i=1\ldots,\ell
\end{array}\right. \quad \Rightarrow \quad\ParV=\ParVT$$
where $D(0 |\varphi_i;\ParV):= \VIPar{0}{\varphi_i}{\ParV}-\VIPar{0}{\varphi_i}{\ParVT}$
\end{itemize}
\item[\textbf{[C3]}] The function $U_n(\varphi;\cdot)$ is continuous for $P_{\ParVT}-$a.e. $\varphi$.
\item[\textbf{[C4]}] For all $\varphi \in \cSp$, $ V^{}\left( 0|\varphi; \ParV \right)$ is continuously differentiable in $\ParV$ and for all $j=1,\ldots,p$
\begin{eqnarray*}
\hspace*{-1cm}\Esp\left( \max_{\ParV \in \SpPar}\left(\left|  {\frac{\partial V^{}}{\partial \theta_{j}}}  \left( 0|\Phi; \ParV \right) \right| e^{- V^{}\left( 0|\Phi; \ParV \right)} \right)^2\right)&<&+\infty .\\
\end{eqnarray*}

\end{itemize}

\begin{theorem} \label{thm-cons}
Under the assumptions \textbf{[Mod-E]} and \textbf{[C]}, for $P_{\ParVT}-$almost every $\varphi$,  the maximum pseudo-likelihood estimate
$\Estn{\ParV}(\varphi)$ converges towards $\ParVT$
as $n$ tends to infinity.
\end{theorem}

\subsection{Asymptotic normality of the MPLE}

For establishing the asymptotic normality of the MPLE we need to assume the four additional following assumptions: 
\begin{itemize}
\item[\textbf{[N1]}] For all $\varphi \in \cSp$,
$\VIPar{0}{\varphi}{\ParV}$ is differentiable in $\ParV=\ParVT$. For all $k=1,\ldots,3$ and for all $\lambda_1,\ldots,\lambda_{k}$, $k$ positive integers such that $\sum_{i=1}^{k} \lambda_i=3$ and for  $\Delta \Subset \RR[2]$
$$
\Esp \left( 
\int_{\Delta^{k}} \prod_{i=1}^{k} 
\left|\dVIPar{j}{0^{M}}{\Phi}{\ParVT}\right|^{\lambda_i}
e^{-\VIPar{\{x_1,\ldots,x_{k}\}}{\Phi}{\ParVT}} dx_1\ldots dx_{k}
\right) <+\infty.
$$

\item[\textbf{[N2]}] There exists a neighbourhood $\Vois(\ParVT)$ of $\ParVT$ such that for all $\varphi\in \cSp$, $\VIPar{0}{\varphi}{\ParV}$ is twice continuously differentiable in $\ParV\in \Vois$ and,  for all $j,k=1,\ldots,p$ and $\ParV\in \Vois(\ParVT)$,
$$\Esp\left( \left| \ddVIPar{j}{k}{0}{\Phi}{\ParV}\right| e^{-\VIPar{0}{\Phi}{\ParV} } 
\right)<+\infty, \quad  \Esp\left( \left| \ddVIPar{j}{k}{0}{\Phi}{\ParV}\right| e^{-\VIPar{0}{\Phi}{\ParVT} } 
\right)<+\infty, $$
and 
$$\Esp\left( \left( \left| \dVIPar{j}{0}{\Phi}{\ParV}\right| e^{-\VIPar{0}{\Phi}{\ParV}} \right)^2
\right)<+\infty.$$
\item[\textbf{[N3]}] There exists 
$A_1,\ldots,A_{\ell}$, $\ell\geq p$ events in $\cSp$ such that:
\begin{itemize}
\item the $\ell$ events $A_i$ are disjoint and satisfy $P_{\ParVT}(A_i) >0$ 
\item for all $\left(\varphi_1,\ldots,\varphi_\ell\right) \in A_1 \times\cdots \times A_{\ell}$ the $(\ell, p)$ matrix with entries $\dVIPar{j}{0}{\varphi_i}{\ParVT}$ is injective.
\end{itemize}
\item[\textbf{[N4]}] There exists $A_0,\ldots,A_{\ell}$, $\ell\geq p$ disjoint sub-events of $\cSptt:=\left\{\varphi \in \cSp: \varphi_{\Dom[i]}=\emptyset, 1\leq |i| \leq 2 \right\}$ such that 
\begin{itemize}
\item for $j=0,\ldots,\ell$, $P_{\ParVT}(A_j)>0$.
\item for all $\left(\varphi_0,\ldots,\varphi_{\ell} \right)\in A_0 \times \cdots \times A_{\ell}$ the $(\ell, p)$ matrix with entries $\left( \ensuremath{\Vect{LPL}_{\Ltt}^{(1)} \left( \varphi_i ; \ParVT  \right)}  \right)_j  - \left(\ensuremath{\Vect{LPL}_{\Ltt}^{(1)} \left( \varphi_0 ; \ParVT  \right)} \right)_j$ is injective,
 with $\Ltt:=\cup_{i \in \BDt[0]{1} }$.
\end{itemize}
\end{itemize}

The assumptions \textbf{[N3]} and \textbf{[N4]} will ensure (see Section~\ref{sec-proofs} for more details) that the matrices $\ddU[\ParVT]$ and $\Mat{\Sigma} (\widetilde{D},\ParVT )$ respectively defined by
\begin{equation}\label{def-U2}
\left(\Mat{U}^{(2)}(\ParVT)\right)_{j,k} := \Esp \left(  {\frac{\partial V^{}}{\partial \theta_{j}}}  \left( 0|\Phi; \ParVT \right)  {\frac{\partial V^{}}{\partial \theta_{k}}}  \left( 0|\Phi; \ParVT \right)
e^{- V^{}\left( 0|\Phi; \ParVT \right)}\right)
\end{equation}
and
\begin{equation} \label{eq-defSig}
\Mat{\Sigma} ( \Dt,\ParVT ) = \Dt^{-2} \sum_{i \in \BDt[0]{1}} \Esp \left(
\ensuremath{\Vect{LPL}_{\Dom[0]}^{(1)} \left( \Phi ; \ParVT  \right)} \tr{\ensuremath{\Vect{LPL}_{\Dom[i]}^{(1)} \left( \Phi ; \ParVT  \right)} }
\right),
\end{equation}
are definite positive.


Observe that,  when the energy function is linear, the expressions of the assumptions \textbf{[N1]} and \textbf{[N2]} are clearly simpler (see \cite{A-BilCoeDro08}) and that \textbf{[C2]} and \textbf{[N3]} are similar.

\begin{theorem} \label{thm-norm}
Under the assumptions \textbf{[Mod]}, \textbf{[C]}, \textbf{[N1]}, \textbf{[N2]} and \textbf{[N3]}, we have the following convergence in distribution as $n \to +\infty$
\begin{equation} \label{convLoiMPLE1}
{}|\Lambda_n|^{1/2}\; \ddU[\ParVT] \; \left( \Estn{\ParV}(\Phi) - \ParVT  \right) \rightarrow  \mathcal{N} \left( 0 , \Mat{\Sigma}(\Dt,\ParVT) \right),
\end{equation}
where $\Mat{\Sigma}(\Dt,\ParVT)$ is defined by~\eqref{eq-defSig}. In addition under the assumption \textbf{[N4]}
\begin{equation} \label{convLoiMPLE2}
{}|\Lambda_n|^{1/2} \; \ensuremath{\Estn{\Mat{\Sigma}}(\Phi;\Estn{\ParV}(\Phi))}^{-1/2}
\; \ensuremath{\Mat{U}_n^{(2)}(\Phi;\Estn{\ParV}(\Phi))}
 \; \left( \Estn{\ParV}(\Phi) - \ParVT  \right) \rightarrow 
\mathcal{N} \left( 0 , \Mat{I}_{p} \right),
\end{equation}
where for some $\ParV$ and any configuration $\varphi$, the matrix \ensuremath{\Estn{\Mat{\Sigma}}(\varphi;\ParV)}  is defined by
\begin{equation} \label{def-EstSign}
\ensuremath{\Estn{\Mat{\Sigma}}(\varphi;\ParV)}  = |\Lambda_n|^{-1} \sum_{i \in I_n}\sum_{j\in  \BDt[i]{1}\cap I_n}   \ensuremath{\Vect{LPL}_{\Delta_i}^{(1)} \left( \varphi ; \ParV  \right)} \tr{\ensuremath{\Vect{LPL}_{\Delta_j}^{(1)} \left( \varphi ; \ParV  \right)}}.
\end{equation}
\end{theorem}

\noindent In the following the assumption \textbf{[N]} will stand for the assumptions \textbf{[N1]}, \textbf{[N2]}, \textbf{[N3]} and \textbf{[N4]}. 



\section{Applications to the LJ-type model} \label{sec-LJ}

This section focuses on the LJ-type model presented in Section~\ref{sec-LJtype} and aims at proving the following result.

\begin{proposition} \label{prop-LJ} ${ }$\\
$(i)$ Theorem~\ref{thm-cons} holds for the LJ-type model (with $D\in ]0,+\infty]$), that is for the Lennard-Jones and the finite-range Lennard-Jones model.\\
$(ii)$ Theorem~\ref{thm-norm} holds only for the finite-range Lennard-Jones model.
\end{proposition}


The proof of Proposition~\ref{prop-LJ} consists in verifying Assumptions~\textbf{[C]} for the LJ-type model and~\textbf{[N]} only for the finite range Lennard-Jones model. In the following, we will deal with two types of assumptions:
\begin{itemize}
\item Integrabilility type assumptions, i.e. Assumptions \textbf{[C1]}, \textbf{[C4]}, \textbf{[N1]} and \textbf{[N2]}. 
\item Identifiability type assumptions, i.e. Assumptions \textbf{[C2]}, \textbf{[N3]} and \textbf{[N4]}.
\end{itemize}
Note that \textbf{[C3]} is obvious since $g^{LJ}(r,\cdot)$ is continuous. For the integrability type assumptions, the following Lemma will be widely used.

\begin{lemma} \label{lem-integ} Let $\Phi$ be a stationary pairwise interaction Gibbs point process assumed to be superstable, regular and lower regular. For $i=1,2$, define $ H_{i}\left(x|\varphi \right)= \sum_{y\in \varphi} g_i(||x-y||)$ with $g_i$ a continuous function. Assume that there exists $\varepsilon>0$ such that there exists a positive and decreasing function $g(\cdot)$ such that $g_{\varepsilon}(r):=g_2(r)- \varepsilon | g_1(r)|\geq -g(r)$ for all $r>0$ and $\int_0^{+\infty} rg(r)dr<+\infty$. Then for all $k\geq 0$,
$$
\Esp\left( 
\left|  H_{1}\left(0|\Phi \right) \right|^k e^{- H_{2}\left(0|\Phi \right)}
\right) <+\infty.
$$
\end{lemma}
\begin{proof}
For all finite configuration $\varphi$
\begin{eqnarray*}
\left|  H_{1}\left(0|\varphi \right) \right|^k e^{- H_{2}\left(0|\varphi \right)} &=&\left|  H_{1}\left(0|\varphi \right) \right|^k \; e^{-\varepsilon | H_{1}\left(0|\varphi \right) |} \; e^{-\left(  H_{2}\left(0|\varphi \right)-\varepsilon  H_{1}\left(0|\varphi \right)\right)} \\
&\leq & c(\varepsilon,k) e^{-\left(  H_{2}\left(0|\varphi \right)-\varepsilon  H_{1}\left(0|\varphi \right)\right)}, \quad \mbox{ with } c(\varepsilon,k)=\left( \frac{k}{\varepsilon e}\right)^k \\
&\leq & c(\varepsilon,k) e^{-  H_{\varepsilon}\left(0|\varphi \right) },
\end{eqnarray*}
where 
$$
 H_{\varepsilon}\left(0|\varphi \right) := \sum_{x\in \varphi} g_{\varepsilon}(||x||).$$
Now, the assumptions ensure that $g_{\varepsilon}$ is lower regular in the Ruelle sense. We may now apply the same argument as in Lemma~3 of~\cite{Mase95} to prove the integrability of the random variable $e^{-  H_{\varepsilon}\left(0|\Phi \right) }$.
\end{proof}

Before verifying the different assumptions, let us denote by  
$$
{\theta_{i}^{\inf{}}}:=\inf_{\ParV\in \SpPar} \theta_i, \quad {\theta_{i}^{\sup}}:=\sup_{\ParV\in \SpPar} \theta_i,\quad {\theta_{}^{\inf{}}}:=\min({\theta_{2}^{\inf{}}},{\theta_{3}^{\inf{}}}) \quad \mbox{ and } \quad {\theta_{}^{\sup}}:=\max({\theta_{2}^{\sup}},{\theta_{3}^{\sup}}).$$
Since $\SpPar$ is a compact set of $\RR[] \times (]0,+\infty[)^2$, then ${\theta_{}^{\inf{}}}>0$ and ${\theta_{}^{\sup}}<+\infty$.

\subsection{Assumptions $\textbf{[C]}$}

\subsubsection{Assumption \textbf{[C1]} }

The first part is a direct application of Lemma~\ref{lem-integ}. For the second part, one has to prove that for all $\ParV \in \SpPar$
$$
\Esp \left(  | H^{LJ}\left( 0|\Phi; \ParV \right) | e^{- H^{LJ}\left( 0|\Phi; \ParVT \right) }\right) <+\infty
$$
Let $g_{\varepsilon}(r)=g^{LJ}(r;\ParVT) - \varepsilon | g^{LJ}(r;\ParV)|$. We have
$$
g_{\varepsilon}(r) :=\left\{ \begin{array}{ll}
4\theta_2^\star \left( \frac{(\theta_3^\star)^{12}-\varepsilon\frac{\theta_2}{\theta_2^\star}\theta_3^{12}}{r^{12}} -\frac{(\theta_3^\star)^{6}-\varepsilon\frac{\theta_2}{\theta_2^\star}\theta_3^{6}}{r^{6}} \right) & \mbox{ if } r\leq \theta_3 \\
 4\theta_2^\star \left( \frac{(\theta_3^\star)^{12}+\varepsilon\frac{\theta_2}{\theta_2^\star}\theta_3^{12}}{r^{12}} -\frac{(\theta_3^\star)^{6}+\varepsilon\frac{\theta_2}{\theta_2^\star}\theta_3^{6}}{r^{6}} \right)& \mbox{ if } r\geq \theta_3
\end{array}\right.
$$
which satisfies the assumptions of Lemma~\ref{lem-integ} as soon as $\varepsilon<\left(\frac{\theta_3^\star}{\theta_3}\right)^{12} \frac{\theta_2^\star}{\theta_2}$, that is, as soon as $\varepsilon<\left(\frac{{\theta_{}^{\inf{}}}}{{\theta_{}^{\sup}}}\right)^{13}$.

\subsubsection{Assumption \textbf{[C2]} }

Let us denote for $n\geq 1$, $C_n = \mathcal{B}(0,n)\setminus \mathcal{B}(0,n-1)$ and define for $m,n\geq 1$ the following configuration sets
\begin{eqnarray*}
U_{m,n} &=&\left\{
\varphi \in \Omega : |\varphi_{C_n}| \leq m |C_n| \right\} \\
U_m &=& \cap_{n\geq 1} U_{m,n}.
\end{eqnarray*}
In order to prove \textbf{[C2]}, we need the following Lemma.
\begin{lemma} \label{lem-Um}
Let $R \in {\RR[]}^+$, $\ParV \in \SpPar$ and $\varphi\in U_m$, let us denote by 
$$Z(\varphi,R;\ParV):=\sum_{x\in \varphi_{\mathcal{B}(0,R)^c}} g^{LJ}(||x||;\ParV),
$$
then for all $\delta>0$ there exists $R_0$ such that for all $R\geq R_0$, $|Z(\varphi,R;\ParV)|\leq \delta$.
\end{lemma}

\begin{proof}
\begin{eqnarray*}
Z(\varphi,R;\ParV)=\big| \sum_{x\in \varphi_{\mathcal{B}(0,R)^c}} g^{LJ}(||x||;\ParV) \big| &\leq & \sum_{n\geq \lceil R \rceil} \sum_{x \in \varphi_{C_n}}\left| g^{LJ}(||x||;\ParV) \right| \\
&\leq & \sum_{n\geq \lceil R \rceil} |\varphi_{C_n}|\times \sup_{x\in C_n}\left|g^{LJ}(||x||;\ParVT) \right|.
\end{eqnarray*}
There exists a constant $k=k(R)$ such that for all $n\geq  \lceil R \rceil$, $\sup_{x\in C_n}\left|g^{LJ}(||x||;\ParVT)\right| \leq k n^{-6}$. Therefore,
$$\big| \sum_{x\in \varphi_{\mathcal{B}(0,R)^c}} g^{LJ}(||x||;\ParV) \big| \leq  k m \sum_{n\geq  \lceil R \rceil} |C_n|\times n^{-6} = \mathcal{O} \left(\sum_{n\geq  \lceil R \rceil}n^{-5} \right),
$$
which leads to the result since the previous series is convergent.
\end{proof}

Let $\ParV \in \SpPar\setminus\ParVT$ and consider the following configuration sets defined for $k\geq 1$ and for $\eta$ small enough by
\begin{eqnarray}
A_0 &=& \left\{ \varphi \in \Omega : | \varphi\cap \mathcal{B}(0,D) |=0 \right\} \label{set1} \\
A_k(\eta) &=& \left\{ \varphi\in \Omega: |\varphi \cap \mathcal{B}(0,D) |= |\varphi\cap \mathcal{B}( (0,D k^{-1/12}),\eta)|=1
\right\}, \label{set2}
\end{eqnarray}
where $D$ is any positive real for the Lennard-Jones model and corresponds to the range of the function $g^{LJ}(\cdot)$ for the finite range Lennard-Jones model.
There exists $m\geq 1$ such that for all $\eta>0$ and for $k=2,4$
$$
P_{\ParVT} \left( A_0 \cap U_m\right)>0
\quad \mbox{ and } \quad
P_{\ParVT} \left( A_k(\eta) \cap U_m\right)>0.
$$
Now, let $\varphi_0 \in A_0\cap U_m, \varphi_2 \in A_2(\eta)\cap U_m$ and $\varphi_4 \in A_4(\eta)\cap U_m$. 
First,
$$
D(0|\varphi_0;\ParV) = \theta_1-\theta_1^\star + Z(\varphi_0,D;\ParV)-Z(\varphi_0,D;\ParVT)=0.$$
For the Lennard-Jones model, according to Lemma~\ref{lem-Um} one has, for $D$ large enough,
$$
\left| Z(\varphi_0,D;\ParV)- Z(\varphi_0,D;\ParVT)\right|\leq \frac12 \left|\theta_1 - \theta_1^\star\right|.
$$
Hence for $\eta$ small enough, and for both models
\begin{eqnarray*}
0&=&\left|D(0|\varphi_0;\ParV) \right| \\
&\geq& |\theta_1-\theta_1^\star| - \left| Z(\varphi_0,D;\ParV)- Z(\varphi_0,D;\ParVT)\right| \\
&\geq& \frac12 |\theta_1-\theta_1^\star|,
\end{eqnarray*}
which leads to $\theta_1=\theta_1^\star$. Moreover,
\begin{eqnarray*}
D(0|\varphi_2;\ParV) &=& 4\theta_2 \left( 2\left( \frac{\theta_3}{D}\right)^{12} -\sqrt{2} \left(\frac{\theta_3}{D}\right)^{6}\right) -4\theta_2^\star \left( 2\left( \frac{\theta_3^\star}{D}\right)^{12} -\sqrt{2} \left(\frac{\theta_3^\star}{D}\right)^{6}\right) \\
&&+f_2(\varphi_2)+ Z\left(\varphi_2,D;\ParV\right)-Z\left(\varphi_2,D;\ParVT\right) \\
D(0|\varphi_4;\ParV) &=& 4\theta_2 \left( 4\left( \frac{\theta_3}{D}\right)^{12} -2 \left(\frac{\theta_3}{D}\right)^{6}\right) -4\theta_2^\star \left( 4\left( \frac{\theta_3^\star}{D}\right)^{12} -2 \left(\frac{\theta_3^\star}{D}\right)^{6}\right) \\
&&+f_4(\varphi_4)+ Z\left(\varphi_4,D;\ParV\right)-Z\left(\varphi_4,D;\ParVT\right), 
\end{eqnarray*}
where for any $\varphi_k \in A_k(\eta)$ ($k=2,4$), there exists a positive function $\widetilde{f}_k(\eta)$ converging towards zero as $\eta\to 0$ such that $|f_k(\varphi_k)|$ is bounded by $\widetilde{f}_k(\eta)$. Now, we have
\begin{eqnarray*}
2 D(0|\varphi_2;\ParV)- D(0|\varphi_4;\ParV) &=& \frac{4(2-2\sqrt{2})}{D^6} \left( \theta_2 \theta_3^6 - \theta_2^\star {\theta_3^\star}^6 \right) + 2f(\varphi_2)-f_4(\varphi_4)+ Z^\prime(\varphi_2,\varphi_4,D;\ParV,\ParVT)\\
&=&0 
\end{eqnarray*}
with 
$$Z^\prime(\varphi_2,\varphi_4,D;\ParV,\ParVT):=2 \left(Z(\varphi_2,D;\ParV)-Z(\varphi_2,D;\ParVT)\right) -\left(Z(\varphi_4,D;\ParV)-Z(\varphi_4,D;\ParVT)\right).
$$
For $\eta$ small enough, we have, for any $\varphi_k \in A_k(\eta)$ ($k=2,4$), 
$$| 2f(\varphi_2)-f_4(\varphi_4)|\leq 2\widetilde{f}_2(\eta)+\widetilde{f}_4(\eta)\leq \frac14 \left| \frac{4(2-2\sqrt{2})}{D^6} \right| |\theta_2 \theta_3^6 - \theta_2^\star {\theta_3^\star}^6|.$$
For the finite range Lennard-Jones model, $Z^\prime(\varphi_2,\varphi_4,D;\ParV,\ParVT)=0$. For the  Lennard-Jones model, according to Lemma~\ref{lem-Um}, one has for $D$ large enough
$$
\left| Z^\prime(\varphi_2,\varphi_4,D;\ParV,\ParVT)\right|\leq \frac14 \left| \frac{4(2-2\sqrt{2})}{D^6} \right| |\theta_2 \theta_3^6 - \theta_2^\star {\theta_3^\star}^6|.
$$
Hence for $\eta$ small enough, and for both models
\begin{eqnarray*}
0&=&\left|\frac{4(2-2\sqrt{2})}{D^6} \left( \theta_2 \theta_3^6 - \theta_2^\star {\theta_3^\star}^6 \right) + 2f(\varphi_2)-f_4(\varphi_4)+ Z^\prime(\varphi_2,\varphi_4,D;\ParV,\ParVT) \right| \\
&\geq& \left|\frac{4(2-2\sqrt{2})}{D^6} \right| |\theta_2 \theta_3^6 - \theta_2^\star {\theta_3^\star}^6| - | 2f(\varphi_2)-f_4(\varphi_4)|-\left| Z^\prime(\varphi_2,\varphi_4,D;\ParV,\ParVT)\right| \\
&\geq& \frac12\left|\frac{4(2-2\sqrt{2})}{D^6} \right| |\theta_2 \theta_3^6 - \theta_2^\star {\theta_3^\star}^6| 
\end{eqnarray*}
leading to $\theta_2 \theta_3^6 = \theta_2^\star {\theta_3^\star}^6$. By considering the combination $\sqrt{2} D(0|\varphi_2;\ParV)-D(0|\varphi_4;\ParV)$ and using similar arguments as previously, one obtains: $\theta_2 \theta_3^{12} = \theta_2^\star {\theta_3^\star}^{12}$. By computing the ratio of the two last equations, one obtains $\theta_3=\theta_3^\star$ and then $\theta_2=\theta_2^\star$.

\subsubsection{Assumption \textbf{[C4]} }

For all $\varphi \in \cSp$ and for any $\ParV \in \SpPar$,
$ V^{LJ}\left( 0|\varphi; \ParV \right)$ is clearly differentiable in $\ParV$. First, note that \textbf{[C4]} is trivial for $j=1$. For $j=2,3$, let us define:
\begin{eqnarray*}
X_j(\varphi;\ParV)&:= \left| {\frac{\partial V^{LJ}}{\partial \theta_{j}}}  \left( 0|\varphi; \ParV \right) \right|e^{ -  V^{LJ}\left( 0|\varphi; \ParV \right)}. \\
\end{eqnarray*} 
Our aim will be to prove that for $j=2,3$ and for all $k>0$
\begin{equation}\label{integ-grad}
\Esp \left( \max_{\ParV\in \SpPar} X_j(\Phi;\ParV)^k \right) <+\infty.
\end{equation}
In particular, the Assumption \textbf{[C4]} corresponds to~(\ref{integ-grad}) with $k=2$. Let us notice that for all $\varphi\in \Omega$ and for all $\ParV \in \SpPar$
$$
 V^{LJ}\left( 0|\varphi; \ParV \right) \geq V^{\inf}(0|\varphi):={\theta_{}^{\inf{}}}+\sum_{x \in \varphi} g^{\inf}(||x||),
$$
with for some $r>0$, $g^{\inf}(r):=4{\theta_{}^{\inf{}}}\left( \frac{\left({\theta_{}^{\inf{}}}\right)^{12}}{r^{12}}-\frac{\left({\theta_{}^{\sup}}\right)^6}{r^6}\right)$. Let us also underline that for $j=2,3$ 
$$
\frac{\partial g^{LJ}}{\partial \theta_j}(r;\ParV) \geq \widetilde{g}^{\inf}_j(r)
\quad \mbox{ with  }
\widetilde{g}^{\inf}_j(r) := \left\{ \begin{array}{ll}
4\left( \frac{\left({\theta_{}^{\inf{}}}\right)^{12}}{r^{12}}-\frac{\left({\theta_{}^{\sup}}\right)^6}{r^6}\right) & \mbox{ if } j=2,\\
4m\left( \frac{12\left({\theta_{}^{\inf{}}}\right)^{11}}{r^{12}}- \frac{6\left({\theta_{}^{\sup}}\right)^5}{r^6}\right) & \mbox{ if } j=3.
\end{array} \right.
$$
Therefore, by defining $\widetilde{V}^{\inf}_j(0|\varphi):=\sum_{x\in \varphi} \widetilde{g}_j^{\inf}(||x||)$, the result~(\ref{integ-grad}) will be ensured by proving
$$
\Esp \left(\widetilde{V}^{\inf}_j(0|\Phi)e^{-V^{\inf}(0|\Phi)} \right) <+\infty.
$$
According to Lemma~\ref{lem-integ}, in order to prove this, let us denote by $g_{j,\varepsilon}(\cdot)$ the function defined for $j=2,3$, for some $\varepsilon>0$ and for $r>0$ by $g_{j,\varepsilon}(r)=\widetilde{g}_j^{\inf}(r)-\varepsilon\left|g^{\inf}(r)\right|$. On the one hand, one has
$$
g_{2,\varepsilon}(r)=\left\{ \begin{array}{ll}
4 \left( \frac{\left({\theta_{}^{\inf{}}}\right)^{13}-\varepsilon \left({\theta_{}^{\inf{}}}\right)^{12}}{r^{12}} - \frac{{\theta_{}^{\inf{}}} \left({\theta_{}^{\sup}}\right)^6-\varepsilon  \left({\theta_{}^{\sup}}\right)^6}{r^6} \right) & \mbox{ if } r \leq \frac{\left({\theta_{}^{\inf{}}}\right)^2}{\theta_{}^{\sup}}, \\
4 \left( \frac{\left({\theta_{}^{\inf{}}}\right)^{13}+\varepsilon \left({\theta_{}^{\inf{}}}\right)^{12}}{r^{12}} - \frac{{\theta_{}^{\inf{}}} \left({\theta_{}^{\sup}}\right)^6+\varepsilon  \left({\theta_{}^{\sup}}\right)^6}{r^6} \right) & \mbox{ if } r \geq \frac{\left({\theta_{}^{\inf{}}}\right)^2}{\theta_{}^{\sup}},
\end{array} \right.
$$
which satisfies the assumptions of Lemma~\ref{lem-integ} as soon as $\varepsilon<{\theta_{}^{\inf{}}}$. On the other hand
$$
g_{3,\varepsilon}(r)=\left\{ \begin{array}{ll}
 4{\theta_{}^{\inf{}}} \left(\frac{\left({\theta_{}^{\inf{}}}\right)^{12}-12\varepsilon \left({\theta_{}^{\inf{}}}\right)^{11}}{r^{12}}- \frac{\left({\theta_{}^{\sup}}\right)^6-6\varepsilon \left({\theta_{}^{\sup}}\right)^5}{r^6} \right)& \mbox{ if } r \leq \left(2 \frac{\left({\theta_{}^{\inf{}}}\right)^{11}}{\left({\theta_{}^{\sup}}\right)^5}\right)^{1/6}\\
 4{\theta_{}^{\inf{}}} \left(\frac{\left({\theta_{}^{\inf{}}}\right)^{12}+12\varepsilon \left({\theta_{}^{\inf{}}}\right)^{11}}{r^{12}}- \frac{\left({\theta_{}^{\sup}}\right)^6+6\varepsilon \left({\theta_{}^{\sup}}\right)^5}{r^6} \right)& \mbox{ if } r \geq \left(2 \frac{\left({\theta_{}^{\inf{}}}\right)^{11}}{\left({\theta_{}^{\sup}}\right)^5}\right)^{1/6},
\end{array} \right.
$$
which satisfies the assumptions of Lemma~\ref{lem-integ} as soon as $\varepsilon<{\theta_{}^{\inf{}}}/12$, which ends the proof.

\subsection{Assumptions \textbf{[N]} }

\subsubsection{Assumption \textbf{[N1]} }

Let us present two auxiliary lemmas.
\begin{lemma} \label{lem-xK} Let $\varphi$ be the realization of a stationary pairwise interaction point process with local energy function defined by
$$
 V^{}\left( x|\varphi; \ParV \right) = \theta_1 +  H^{}\left( x|\varphi; \ParV \right) \quad \mbox{ with } \quad H^{}\left( x|\varphi; \ParV \right)= \sum_{y\in \varphi} g(||y-x||;\ParV).
$$
Let $K<+\infty$ and let $x_1,\ldots,x_K \in \RR \setminus \varphi$, $x_i\neq x_j$ for $i,j=1,\ldots,K$ (where $K<+\infty$), then 
\begin{eqnarray*}
 H^{}\left( \{x_1,\ldots,x_K\}|\varphi; \ParV \right) &=& \sum_{k=1}^K  H^{}\left( x_k|\varphi; \ParV \right) +  H^{}\left(\{x_1,\ldots,x_K\}; \ParV \right) \\
 V^{}\left( \{x_1,\ldots,x_K\}|\varphi; \ParV \right) &=& \sum_{k=1}^K  V^{}\left( x_k|\varphi; \ParV \right) +  H^{}\left(\{x_1,\ldots,x_K\}; \ParV \right)
\end{eqnarray*}
\end{lemma}

This result comes from the definition of the local energy.

\begin{lemma} \label{lem-aux2}
Using the same notation and under the same assumptions of Lemma~\ref{lem-xK}, assume that there exists $g_{min}$ such that for all $r>0$ and any $\ParV \in \SpPar$, $g(r;\ParV)\geq g_{min}$, then
$$
e^{- V^{}\left( \{x_1,\ldots,x_K\}|\varphi; \ParV \right)} \leq c_K \prod_{k=1}^K e^{- V^{}\left( x_k|\varphi; \ParV \right)} \quad \mbox{ with } c_K=e^{ -\frac{K(K-1)}{2} g_{min}}
$$
\end{lemma}

\begin{proof}
The proof is immediate since
$$
 H^{}\left(\{x_1,\ldots,x_K\}; \ParV \right) = \sum_{i<j} g(||x_i-x_j||;\ParV) \geq \frac{K(K-1)}2 g_{min}.
$$
\end{proof}

Let $k=1,\ldots,3$ and let $\lambda_1,\ldots,\lambda_{k}$, $k$ positive integers such that $\sum_{i=1}^{k} \lambda_i=3$ and define the random variable
$$
A(\Phi) := \int_{\Delta^{k}} \prod_{i=1}^{k} 
\left|\dVIPar{j}{x_i^{}}{\Phi}{\ParVT}\right|^{\lambda_i}
e^{-\VIPar{\{x_1,\ldots,x_{k}\}}{\Phi}{\ParVT}}dx_i.
$$
From Lemma~\ref{lem-aux2}, we have
\begin{eqnarray*}
\Esp\left( A(\Phi)\right) &\leq & \Esp\left( c_{k}  \int_{\Delta^{k}}\prod_{i=1}^{k} 
\left|\dVIPar{j}{x_i^{}}{\Phi}{\ParVT}\right|^{\lambda_i} e^{-\VIPar{x_i}{\Phi}{\ParVT}} dx_i \right) \\
&=&  c_{k}  \int_{\Delta^{k}}\Esp\left(\prod_{i=1}^{k} 
\left|\dVIPar{j}{x_i^{}}{\Phi}{\ParVT}\right|^{\lambda_i} e^{-\VIPar{x_i}{\Phi}{\ParVT}} \right) dx_1\ldots dx_{k} \\
&\leq & c_{k}    \int_{\Delta^{k}} \prod_{i=1}^{k} \Esp \left(
\left|\dVIPar{j}{x_i^{}}{\Phi}{\ParVT}\right|^{k} e^{-\frac{k}{\lambda_i}\VIPar{x_i}{\Phi}{\ParVT}} 
\right)^{1/k} dx_1\ldots dx_k
\\
&=& c_{k}   \prod_{i=1}^{k} \int_{\Delta}  \Esp \left(
\left|\dVIPar{j}{x_i^{}}{\Phi}{\ParVT}\right|^{k} e^{-\frac{k}{\lambda_i}\VIPar{x_i}{\Phi}{\ParVT}} 
\right)^{1/k} dx_i \\
&=& c_{k}  |\Delta|^k \prod_{i=1}^{k} \Esp \left(
\left|\dVIPar{j}{0}{\Phi}{\ParVT}\right|^{k} e^{-\frac{k}{\lambda_i}\VIPar{0}{\Phi}{\ParVT}} 
\right)^{1/k}
\end{eqnarray*}
by using Hölder's inequality and the stationarity of the process. The result is then a simple consequence of~(\ref{integ-grad}) and Lemma~\ref{lem-integ}.

\subsubsection{Assumption \textbf{[N2]} }

For all $\varphi \in \cSp$, it is clear that for all $\ParV \in \SpPar$,
$\VIPar{0}{\varphi}{\ParV}$ is twice continuously differentiable in $\ParV$. According to Lemma~\ref{lem-integ} and the fact that \textbf{[N1]} is satisfied, it is sufficient to prove that for all $j,k=1,2,3$
$$
\Esp\left( \left|  
\frac{\partial^2 V^{LJ}}{\partial \theta_{j}\partial \theta_{k} }  \left( 0|\Phi; \ParV \right)\right| e^{- V^{LJ}\left( 0|\Phi; \ParV \right) } 
\right)<+\infty.
$$
This is obvious when either $j$ or $k$ equals 1 and when $j=k=2$ (since $\frac{\partial^2 g^{LJ}}{(\partial {\theta}_2)^2}(r;\cdot)=0$). Now, for the other cases, define for $\ParV \in \SpPar$
$g_{j,k,\varepsilon}(r):=g^{LJ}(r;\ParV)-\varepsilon\left|  \frac{\partial^2 g^{LJ}}{\partial \theta_j \partial \theta_k}(r;\ParV) \right|$. We have
$$
g_{2,3,\varepsilon}(r)=g_{3,2,\varepsilon}(r)=\left\{ \begin{array}{ll}
4 \left(
\frac{\theta_2 \theta_3^{12}-12 \varepsilon \theta_3^{11}}{r^{12}} - \frac{\theta_3^6-6\varepsilon\theta_3^5}{r^6}
\right) & \mbox{ if } r\leq 2^{1/6} \\
4 \left(\frac{\theta_2 \theta_3^{12}+12 \varepsilon \theta_3^{11}}{r^{12}} - \frac{\theta_3^6+6\varepsilon\theta_3^5}{r^6}
\right) & \mbox{ otherwise}
\end{array}\right.
$$
which satisfies the assumptions of Lemma~\ref{lem-integ} as soon as $\varepsilon < \frac{\theta_2 \theta_3}{12}$, that is, as soon as $\varepsilon<\frac{\left({\theta_{}^{\inf{}}}\right)^2}{12}$. Finally,
$$
g_{3,3,\varepsilon}(r)=\left\{ \begin{array}{ll}
4\left( 
\frac{\theta_2 \theta_3^{12} - 132 \varepsilon\theta_3^{10}}{r^{12}} - \frac{\theta_2 \theta_3^6 - 30 \varepsilon \theta_3^4}{r^6}
\right) & \mbox{ if } r\leq \left( \frac{132}{30}\right)^{1/6}\theta_3 \\
4\left( 
\frac{\theta_2 \theta_3^{12} + 132 \varepsilon\theta_3^{10}}{r^{12}} - \frac{\theta_2 \theta_3^6 + 30 \varepsilon \theta_3^4}{r^6}
\right) & \mbox{ otherwise}
\end{array}\right.
$$
which satisfies the assumptions of Lemma~\ref{lem-integ} as soon as $\varepsilon < \frac{\theta_2\theta_3^2}{132}$, that is, as soon as $\varepsilon < \frac{\left({\theta_{}^{\inf{}}}\right)^3}{132}$.

\subsubsection{Assumption \textbf{[N3]} }

Let $\Vect{y}=(y_1,y_2,y_3)\in \RR[3]$ and $g(\Vect{y},\varphi):=\tr{\Vect{y}}  \Vect{V}^{(1)}_{LJ}\left( 0|\varphi; \ParVT \right)$. Let $\varphi_0 \in A_0$ and $\varphi_k(\eta) \in A_k(\eta)$ ($k=2,4$) where $A_0$ and $A_k(\eta)$ are defined by~(\ref{set1}) and~(\ref{set2}). Assume $g(\Vect{y},\varphi_k)=0$ for $k=0,2,4$. Since, $g(\Vect{y},\varphi_0) =y_1$, we have $y_1=0$. Now, 
\begin{eqnarray*}
g(\Vect{y},\varphi_2) &=& 4y_2 \left( 2 \left( \frac{\theta_3^\star}D\right)^{12}- \sqrt{2} \left(\frac{\theta_3^\star}D \right)^6\right) + 4y_3 \theta_2^\star \left(  2 \frac{12{\theta_3^\star}^{11}}{D^{12}}- \sqrt{2} \frac{6{\theta_3^\star}^5}{D^6}\right) + f_2(\Vect{y},\varphi_2)\\
g(\Vect{y},\varphi_4) &=&4y_2 \left( 4 \left( \frac{\theta_3^\star}D\right)^{12}- 2 \left(\frac{\theta_3^\star}D \right)^6\right) + 4y_3 \theta_2^\star \left(  4 \frac{12{\theta_3^\star}^{11}}{D^{12}}- 2 \frac{6{\theta_3^\star}^5}{D^6}\right) + f_4(\Vect{y},\varphi_4),
\end{eqnarray*}
where for any $\varphi_k \in A_k(\eta)$ ($k=2,4$), there exists a positive function $\widetilde{f}_k(\Vect{y},\eta)$ converging towards zero as $\eta\to 0$ such that $|f_k(\Vect{y},\varphi_k)|$ is bounded by $\widetilde{f}_k(\Vect{y},\eta)$. Now, we have
$$
2g(\Vect{y},\varphi_2)-g(\Vect{y},\varphi_4) = 4(2-2\sqrt{2}))\frac{{\theta_3^\star}^5}{D^6} \left({\theta_3^\star} y_2  + 6 \theta_2^\star y_3 \right) + 2f_2(\Vect{y},\varphi_2)-f_4(\Vect{y},\varphi_4)=0.
$$
For $\eta$ small enough, we have, for any $\varphi_k \in A_k(\eta)$ ($k=2,4$), 
$$| 2f(\Vect{y},\varphi_2)-f_4(\Vect{y},\varphi_4)|\leq 2|\widetilde{f}_2(\Vect{y},\eta)|+|\widetilde{f}_4(\Vect{y},\eta)|\leq \frac12\left|  4(2-2\sqrt{2})\frac{{\theta_3^\star}^5}{D^6} \left(\theta_3^\star y_2  + 6 \theta_2^\star y_3 \right) \right|.$$
Hence for $\eta$ small enough, 
$$
0= \left|2g(\Vect{y},\varphi_2)-g(\Vect{y},\varphi_4)\right| \geq \frac12 \left|  4(2-2\sqrt{2})\frac{{\theta_3^\star}^5}{D^6} \left(\theta_3^\star y_2  + 6 \theta_2^\star y_3 \right) \right|,
$$
leading to the equation $\theta_3^\star y_2  + 6 \theta_2^\star y_3=0$. By considering the linear combination $\sqrt{2}g(\Vect{y},\varphi_2)-g(\Vect{y},\varphi_4)$, we may obtain the equation $\theta_3^\star y_2  + 12 \theta_2^\star y_3=0$ with similar arguments. Both equations lead to $y_2=y_3=0$.

\subsubsection{Assumption \textbf{[N4]} }

The assumption \textbf{[N4]} may be rewritten for all $k=1,\cdots,\ell$ and for all $\varphi_k \in A_k$ and $\varphi_0\in A_0$:
$$
\left(\forall \Vect{y}\in \RR[3], \tr{\Vect{y}} \left( \ensuremath{\Vect{LPL}_{\Ltt}^{(1)} \left( \varphi_k ; \ParVT  \right)}-\ensuremath{\Vect{LPL}_{\Ltt}^{(1)} \left( \varphi_0 ; \ParVT  \right)}\right)
 =\tr{\Vect{y}}(\Vect{L}(\varphi_k;\ParVT)- \Vect{R}(\varphi_k;\ParVT)) =0 \right)\Longrightarrow \Vect{y}=0.
$$
where for any configuration $\varphi \in \cSptt$ and $\varphi_0 \in A_0$
\begin{eqnarray*}
\Vect{L}(\varphi;\ParVT)&:=&\int_{\Ltt}{   \Vect{V}^{(1)}_{LJ}\left( x|\varphi; \ParVT \right) e^{-  V^{LJ}\left( x|\varphi; \ParVT \right)} dx}-
\int_{\Ltt}{    \Vect{V}^{(1)}_{LJ}\left( x|\varphi_0; \ParVT \right) e^{-  V^{LJ}\left( x|\varphi_0; \ParVT \right)} dx} \\
\Vect{R}(\varphi;\ParVT)&:=&\sum_{x\in \varphi\cap \Ltt}  \Vect{V}^{(1)}_{LJ}\left( x|\varphi\setminus x; \ParVT \right) - \sum_{x\in \varphi_0\cap \Ltt}  \Vect{V}^{(1)}_{LJ}\left( x|\varphi_0\setminus x; \ParVT \right).
\end{eqnarray*}
Concerning this assumption, we choose $\varphi_0 \in A_0 = \left\{ \varphi\in \overline{\Omega} : \varphi_{\Delta_0}=\emptyset\right\}$. Let $\Vect{y} \in \RR[3]$ then
$$
\int_{\Ltt}{ \tr{\Vect{y}} \Vect{V}^{(1)}_{LJ}\left( x|\varphi_0; \ParVT \right) e^{- V^{LJ}\left( x|\varphi_0; \ParVT \right) }dx } = y_1 e^{-\theta_1^\star}
{}\left| \Ltt \right| \quad \mbox{ and } \quad
\sum_{x\in \varphi_0\cap \Ltt} \tr{\Vect{y}} \Vect{V}^{(1)}_{LJ}\left( x|\varphi_0\setminus x; \ParVT \right)=0.
$$
Consider the following configuration set, defined for $\eta,\varepsilon>0$, by
$$
A_2(\eta,\varepsilon)=\left\{
\varphi\in \overline{\Omega}: \varphi_{\Delta_0}=\{z_1,z_2\} \mbox{ where } z_1 \in \mathcal{B}(0,\eta), z_2 \in \mathcal{B}((0,2\eta+\varepsilon),\eta) 
\right\}.
$$
Note that for $z_1 \in \mathcal{B}(0,\eta), z_2 \in \mathcal{B}((0,2\eta+\varepsilon),\eta)$, $\varepsilon\leq ||z_2-z_1||\leq \varepsilon+4\eta$. Let $\varphi_2\in A_2(\eta,\varepsilon)$ and $x\in \Ltt$, then one may prove that for $j=2,3$
\begin{eqnarray*}
 V^{LJ}\left( x|\varphi_2; \ParVT \right)  &=& \theta_1^\star + 2 g^{LJ}(||x||;\ParVT) + f(x,\eta,\varepsilon) \\
 {\frac{\partial V^{LJ}}{\partial \theta_{j}}}  \left( x|\varphi_2; \ParVT \right) &=& 2 \frac{\partial g^{LJ}}{\partial \theta_j}(||x||;\ParVT)+ f_j(x,\eta,\varepsilon)
\end{eqnarray*}
where $f(x,\eta,\varepsilon)$ and $f_j(x,\eta,\varepsilon)$ are such that 
$$\lim_{(\eta,\varepsilon)\to (0,0)} f(x,\eta,\varepsilon)= \lim_{(\eta,\varepsilon)\to (0,0)} f_j(x,\eta,\varepsilon)= 0.$$
On the one hand, one may prove that there exists a function $f_L(\Vect{y},\eta,\varepsilon)$ such that \\ $\lim_{(\eta,\varepsilon)\to (0,0)} f_L(\Vect{y},\eta,\varepsilon)=0$ and such that
$$
\tr{\Vect{y}}\Vect{L}(\varphi_2;\ParVT) = \tr{\Vect{y}} \Vect{I} -y_1e^{-\theta_1^\star}|\Ltt| + f_L(\Vect{y},\eta,\varepsilon)
$$
where 
$$
\Vect{I}:= \int_{\Ltt} \Vect{h}(||x||;\ParVT)e^{-\theta_1^\star-2g^{LJ}(||x||;\ParVT)}dx
\quad \mbox{ and } \quad 
\Vect{h}(r;\ParVT):=\tr{\left(1, 2\frac{\partial g^{LJ}}{\partial \theta_2}(r;\ParVT),2\frac{\partial g^{LJ}}{\partial \theta_3}(r;\ParVT) \right)}.
$$
On the other hand, there exists a function $f_R(\Vect{y},\eta,\varepsilon)$ such that $\lim_{\eta\to 0} f_R(\Vect{y},\eta,\varepsilon)=0$
$$
\tr{\Vect{y}}\Vect{R}(\varphi_2;\ParVT) = 2y_1 +2y_2 4\left( \left(\frac{\theta_3^\star}{\varepsilon} \right)^{12}-\left(\frac{\theta_3^\star}{\varepsilon} \right)^{6}\right) +
2y_3 4\theta_2^\star \left( 
\frac{12{\theta_3^{\star}}^{11}}{\varepsilon^{12}} - \frac{6{\theta_3^\star}^5}{\varepsilon^6}
\right) + f_R(\Vect{y},\eta,\varepsilon).$$
Since 
\begin{eqnarray*}
\varepsilon^{12} \tr{\Vect{y}}\left( \Vect{L}(\varphi_2;\ParVT)-\Vect{R}(\varphi_2;\ParVT)\right) &=&
\varepsilon^{12}\left( \tr{\Vect{y}}\Vect{I} -y_1e^{-\theta_1^\star}|\Ltt|+f_L(\Vect{y},\eta,\varepsilon)-f_R(\Vect{y},\eta,\varepsilon)\right)\\
&&- \varepsilon^6\left( 2y_2 4{\theta_3^\star}^6+2y_34\theta_2^\star 6{\theta_3^\star}5
\right)
+ 2y_2 4{\theta_3^\star}^{12}+2y_3 4\theta_2^\star 12{\theta_3^\star}^{11}.
\end{eqnarray*}
For $\eta$ and $\varepsilon$ chosen small enough, one may prove that
$$
0= \left| \varepsilon^{12} \tr{\Vect{y}}\left( \Vect{L}(\varphi_2;\ParVT)-\Vect{R}(\varphi_2;\ParVT)\right) \right| \geq \frac12 \left|2y_2 4{\theta_3^\star}^{12}+2y_3 4\theta_2^\star 12{\theta_3^\star}^{11} \right|$$
leading to 
\begin{equation}\label{n4-eq1}
2y_2 4{\theta_3^\star}^{12}+2y_3 4\theta_2^\star 12{\theta_3^\star}^{11}=0 \Leftrightarrow \theta_3^\star y_2 + 12 \theta_2^\star y_3=0.
\end{equation}
This means that
$$\tr{\Vect{y}}\Vect{R}(\varphi_2;\ParVT) = 2y_1 -\frac{1}{\varepsilon^6}\left( 2y_2 4{\theta_3^\star}^6+2y_34\theta_2^\star 6{\theta_3^\star}^5\right)+ f_R(\Vect{y},\eta,\varepsilon).$$
With similar arguments, we obtain that 
\begin{equation}\label{n4-eq2}
2y_2 4{\theta_3^\star}^6+2y_34\theta_2^\star 6{\theta_3^\star}^5=0 \Leftrightarrow \theta_3^\star y_2 + 6 \theta_2^\star y_3=0.
\end{equation}
Equations (\ref{n4-eq1}) and (\ref{n4-eq2}) lead to $y_2=y_3=0$.
Now consider the following configuration set defined for some $k\geq 1$ and $\eta>0$
$$
A_k(\eta)=\left\{ \varphi \in \overline{\Omega}:\varphi_{\Delta_0}= | \varphi\cap \mathcal{B}(0,\eta)|=k\right\}
$$
and let $\varphi_k \in A_k(\eta)$. Then, one may prove that there exists a function $\widetilde{f}_L(\Vect{y},\eta)$ such that $\lim_{\eta\to 0}\widetilde{f}_L(\Vect{y},\eta)=0$ and such that 
$$
\tr{\Vect{y}}\left( \Vect{L}(\varphi_k;\ParVT) - \Vect{R}(\varphi_k;\ParVT)\right)= y_1 \int_{\Ltt} e^{-\theta_1^\star}\left(e^{- k g^{LJ}(||x||;\ParVT)}-1\right)dx -ky_1 + \widetilde{f}_L(\Vect{y},\eta) = 0.
$$
Let us denote by $\Lambda_1:=\mathcal{B}(0,\min(\theta_3^\star,D))$ and $\Lambda_2:=\mathcal{B}(0,D)\setminus \Lambda_1$  
Now let us consider two cases. \\
\textbf{Case 1: $\theta_3^\star\leq D$}. First note that for all $x\in \Ltt$, $g^{LJ}(||x||;\ParVT)\geq 0$. Then, for $k$ large enough and for $\eta$ small enough, we have
$$
\left| \frac1k \int_{\Lambda_1} e^{-\theta_1^\star}\left( e^{-kg^{LJ}(||x||;\ParVT)}-1\right)dx\right| \leq  \frac{|\Lambda_1|}k e^{-\theta_1^\star} \leq  \frac14 \quad \mbox{ and } \quad 
\left| \frac1k \widetilde{f}_L(\Vect{y},\eta)\right| \leq  \frac{|y_1|}4.$$
Hence for $k$ large enough and for $\eta$ small enough, we may obtain
\begin{eqnarray*}
0&=&\frac1k\left|\tr{\Vect{y}}\left( \Vect{L}(\varphi_k;\ParVT) - \Vect{R}(\varphi_k;\ParVT)\right)\right| \\
&\geq& |y_1| -\left| y_1  \frac1k \int_{\Lambda_1} e^{-\theta_1^\star}\left( e^{-kg^{LJ}(||x||;\ParVT)}-1\right)dx + \frac1k \widetilde{f}_L(\Vect{y},\eta)\right|\\
&\geq& |y_1|-\frac{|y_1|}4-\frac{|y_1|}4= \frac{|y_1|}2,
\end{eqnarray*}
which leads to $y_1=0$.

\noindent \textbf{Case 2: $\theta_3^\star\geq D$.} First note that for all $x \in \Lambda_2$, 
$$g^{LJ}(||x||;\ParVT) \leq g_m:=g^{LJ}(D;\ParVT)=4\theta_2^\star\left( \left(\frac{\theta_3^\star}D\right)^{12}-\left(\frac{\theta_3^\star}D\right)^{6}\right)<0.$$
On the one hand, for $k$ large enough and for $\eta$ small enough, we may have
$$
\left| \frac1k y_1\int_{\Lambda_1}e^{-\theta_1^\star}\left( e^{-kg^{LJ}(||x||;\ParVT)}-1\right)dx  + \frac1k \widetilde{f}_L(\Vect{y},\eta) -y_1
\right| \leq \frac{|y_1|}2 +|y_1| \leq \frac32 |y_1|.
$$
On the other hand, we have for $k$ large enough
\begin{eqnarray*}
\frac1k \left|y_1 \int_{\Lambda_2}e^{-\theta_1^\star}\left( e^{-kg^{LJ}(||x||;\ParVT)}-1\right)dx\right| &=& \frac{|y_1|}k \int_{\Lambda_2}e^{-\theta_1^\star}\left( e^{-kg^{LJ}(||x||;\ParVT)}-1\right)dx \\
&\geq& \frac{|y_1|}k e^{-\theta_1^\star} |\Lambda_2| \left( e^{-k g_m}-1\right) = |y_1|e^{-\theta_1^\star} \frac{e^{k|g_m|}-1}k\\
&\geq& 2|y_1|.
\end{eqnarray*}
Therefore for $k$ large enough and for $\eta$ small enough, we have
$$
0=\frac1k\left|\tr{\Vect{y}}\left( \Vect{L}(\varphi_k;\ParVT) - \Vect{R}(\varphi_k;\ParVT)\right)\right| \geq 2|y_1| -\frac32 |y_1| = \frac{|y_1|}2,$$
which leads to $y_1=0$.


\section{Annex: proofs of Theorems \ref{thm-cons} and \ref{thm-norm}} \label{sec-proofs}

Let us start by presenting a particular case of the Campbell Theorem combined with the Glötz Theorem that is widely used in our future proofs. 

\begin{corollary} \label{cor-CampGlotz}
Assume that the point process $\Phi$ with probability measure $P$ is stationary. Let $\Lambda \Subset \RR[2]$, $\varphi\in \cSp$ and let $g$ be a  function satisfying $g(x, \varphi)=g(0,\tau_x\varphi)$ for all $x\in \RR[2]$. Define $f(\varphi)=g(0,\varphi)e^{-\VI{0}{\varphi}}$ and assume that $f\in L^1(P)$. Then,
\begin{equation} \label{eq-CamGlotz}
\Esp \left( \sum_{x \in \Phi_\Lambda} g(x,\Phi \setminus x) \right) = 
\Esp \left( \int_\Lambda g(x,\Phi)e^{-\VI{x}{\Phi}}dx\right)
=|\Lambda| \; \Esp\Big(\; g\left(0,\Phi\right) e^{-\VI{0}{\Phi}}\; \Big)
\end{equation}
\end{corollary}

\begin{proof} see Corollary~3 of \cite{A-BilCoeDro08}
\end{proof}

Let us now present a version of an ergodic theorem obtained by \cite{Nguyen79} and widely used in this paper. Let $\Delta_0$ be a fixed bounded domain

\begin{theorem}[\cite{Nguyen79}] \label{thmNguyen}
Let $\{H_G, G\in \mathbf{\mathcal{B}}_b\}$ be a family of random variables, which is covariant, that for all $x \in \RR$, 
$$
H_{\tau_x G}(\tau_x \varphi) = H_G(\varphi), \;\; \mbox{for }a.e.\; \varphi
$$
and additive, that is for every disjoint $G_1,G_2 \in \mathbf{\mathcal{B}}_b$,
$$
H_{G_1\cup G_2} = H_{G_1} + H_{G_2} , \quad a.s.
$$
Let $\mathcal{I}$ be the sub-$\sigma$-algebra of $\cal F$ consisting of translation invariant (with probability 1) sets. Assume there exists a nonnegative and integrable random variable $Y$ such that $|H_G|\leq Y$ a.s. for every convex $G\subset \Delta_0$. Then,
$$
\lim_{n \to +\infty} \frac1{|G_n|} H_{G_n} = \frac1{|\Delta_0|}E (H_{\Delta_0}|\mathbf{\mathcal{I}}) , \quad a.s.
$$
for each regular sequence $G_n \to \RR$.
\end{theorem}


\subsection{Proof of Theorem~\ref{thm-cons}}

Due to the decomposition of stationary measures as a mixture of ergodic measures (see \cite{B-Pre76}), one only needs to prove Theorem~\ref{thm-cons} by assuming that $P_{\ParVT}$ is ergodic. From now on, $P_{\ParVT}$ is assumed to be ergodic. The tool used to obtain the almost sure convergence is a convergence theorem for minimum contrast estimators established by \cite{Guyon92}. 

\noindent We proceed in three stages.

\noindent \textit{Step~1. Convergence of \ensuremath{U_n(\Phi;\ParV)}}.


Decompose $ \ensuremath{U_n(\varphi;\ParV)}= \frac1{|\Lambda_n|} \left( H_{1,\Lambda_n}(\varphi) +  H_{2,\Lambda_n}(\varphi) \right)$ with
$$
H_{1,\Lambda_n}(\varphi)= \ism[\Lambda_n]{ e^{ -\VIPar{x}{\varphi}{\ParV} }} 
\quad \mbox{ and } \quad 
H_{2,\Lambda_n}(\varphi) = \sum_{x \in \Phi_{\Lambda_n}} \VIPar{x}{\varphi\setminus x}{\ParV}.
$$
Under the assumption \textbf{[C1]}, one can apply Theorem~\ref{thmNguyen} (\cite{Nguyen79}) to the process $H_{1,\Lambda_n}$. And from Corollary~\ref{cor-CampGlotz}, we obtain $P_{\ParVT}-$almost surely as $n \to +\infty$
\begin{equation} \label{H1Lambda}
\frac1{|\Lambda_n|} H_{1,\Lambda_n}(\Phi)
\rightarrow \Esp \Big( e^{- \VIPar{ 0 }{\Phi}{\ParV} } \Big). 
\end{equation}
Now, let $G \subset \Delta_0$, we clearly have
$$
| H_{2,G}(\varphi) |\leq \sum_{x \in \varphi_G} | \VIPar{x}{\varphi\setminus x}{\ParV} | \leq
\sum_{x \in \varphi_{\Delta_0}} | \VIPar{x}{\varphi\setminus x}{\ParV} |. $$
Under the assumption \textbf{[Mod]} and from Corollary~\ref{cor-CampGlotz}, we have 
$$\Esp \left( \sum_{x \in \Phi_{\Delta_0}} | \VIPar{x}{\Phi\setminus x}{\ParV} | \right) = |\Delta_0| \Esp \left( | \VIPar{0}{\Phi}{\ParV} | e^{ - \VIPar{0}{\Phi}{\ParVT} }\right) <+\infty
$$
This means that for all $G \subset \Delta_0$, there exists a random variable $Y \in L^1(P_{\ParVT})$ such that $|H_{2,G}(\Phi)|\leq Y$. Thus, under the assumption \textbf{[C1]} and from Theorem~\ref{thmNguyen} (\cite{Nguyen79}) and from Corollary~\ref{cor-CampGlotz}, we have $P_{\ParVT}-$almost surely
\begin{equation} \label{H2Lambda}
 \frac1{|\Lambda_n|} H_{2,\Lambda_n}(\Phi) \rightarrow  
\frac1{|\Delta_0|} \Esp  \Big( \sum_{x \in \Phi_{\Delta_0}} \VIPar{x}{\Phi\setminus x}{\ParV} \Big) = 
\Esp  \left(\VIPar{0}{\Phi}{\ParV} e^{- \VIPar{0}{\Phi}{\ParVT}} \right).
\end{equation}
We have the result by combining (\ref{H1Lambda}) and (\ref{H2Lambda}): $P_{\ParVT}-$almost surely
\begin{equation} \label{defU}
\ensuremath{U_n(\Phi;\ParV)} \rightarrow U(\ParV) =  \Esp \Big(  e^{
- \VIPar{0}{\Phi}{\ParV} } \; + \;
\VIPar{0}{\Phi}{\ParV} e^{ - \VIPar{0}{\Phi}{\ParVT}}
\Big)
\end{equation}

\noindent \textit{Step~2. \ensuremath{U_n(\cdot;\ParV)} is a contrast function} 

Recall that \ensuremath{U_n(\cdot;\ParV)} is a contrast function if there exists a function $K(\cdot,\ParVT)$ (i.e. nonnegative function equal to zero if and only if $\ParV=\ParVT$) such that $P_{\ParVT}-$almost surely $\ensuremath{U_n(\Phi;\ParV)}- \ensuremath{U_n(\Phi;\ParVT)} \rightarrow K(\ParV,\ParVT)$. From Step 1, we have
\begin{equation} \label{contrasteK}
K(\ParV, \ParVT)\!\!=\!\! \Esp \Big(
e^{-\VIPar{0}{\Phi}{\ParVT}} \Big( 
e^{\VIPar{0}{\Phi}{\ParV}-\VIPar{0}{\Phi}{\ParVT}} - \Big(1+\VIPar{0}{\Phi}{\ParV}-\VIPar{0}{\Phi}{\ParVT} \Big)
\Big) \Big). 
\end{equation}
Since the function $t \mapsto e^t -(1+t)$ is nonnegative and is equal to zero if and only if $t=0$, $K(\ParV,\ParVT)\geq 0$ and 
\begin{eqnarray*}
K(\ParV,\ParVT)=0 &\Leftrightarrow& e^{\VIPar{0}{\varphi}{\ParV}-\VIPar{0}{\varphi}{\ParVT}} - \Big(1+\VIPar{0}{\varphi}{\ParV}-\VIPar{0}{\varphi}{\ParVT} \Big)=0\\
&\Leftrightarrow& D\left(0|\varphi;\ParV \right):=\VIPar{0}{\varphi}{\ParV}-\VIPar{0}{\varphi}{\ParVT}=0
\end{eqnarray*}
 for $P_{\ParVT}-\mbox{a.e. } \varphi$. Let us consider the $\ell$ events $A_j$ ($j=1,\ldots,\ell$) defined in Assumption \textbf{[C2]}. The previous equation is at least true for $\varphi_j \in A_j$, which leads under Assumption \textbf{[C2]} to $\ParV=\ParVT$. Therefore, $K(\ParV,\ParVT)=0 \Rightarrow \ParV=\ParVT$. The converse is trivial.

Before ending this step, note that the assumption \textbf{[C3]} asserts that for any $\varphi$, \ensuremath{U_n(\varphi;\cdot)} and $K(\cdot,\ParVT)$ are continuous functions.

\noindent\textit{Step~3. Modulus of continuity.}

The modulus of continuity of the contrast process defined for all $\varphi\in \cSp$ and all $\eta>0$ by 
$$
W_n(\varphi,\eta) = \sup \left\{ 
\Big|\ensuremath{U_n(\varphi;\ParV)} - \ensuremath{U_n(\varphi;\ParV^\prime)} \Big|: \ParV,\ParV^\prime \in \SpPar, || \ParV - \ParV^\prime || \leq \eta
\right\}
$$
is such that there exists a sequence $(\varepsilon_k)_{k \geq 1}$, with $\varepsilon_k \to 0$
as $k \to +\infty$ such that for all $k \geq 1$
\begin{equation} \label{modCont}
P \left( \limsup_{n \to +\infty}  \left( 
W_n \left(\Phi,\frac1k\right) \geq \varepsilon_k
\right)\right) = 0.
\end{equation}
Let us start to write $W_n\left(\varphi,\frac1k \right) \leq W_{1,n}\left(\varphi,\frac1k \right)+W_{2,n}\left(\varphi,\frac1k \right)$ with
\begin{eqnarray*}
W_{1,n}\left(\varphi,\frac1k \right) &:=& \sup \left\{ 
W_{1,\Lambda_n}^\prime(\varphi;\ParV,\ParV^\prime)
: \ParV,\ParV^\prime \in \SpPar, || \ParV - \ParV^\prime || \leq \frac1k
\right\} \\
W_{2,n}\left(\varphi,\frac1k \right) &:=& \sup \left\{ 
W_{2,\Lambda_n}^\prime(\varphi;\ParV,\ParV^\prime) : \ParV,\ParV^\prime \in \SpPar, || \ParV - \ParV^\prime || \leq \frac1k
\right\}.
\end{eqnarray*}
and
\begin{eqnarray*}
W_{1,\Lambda_n}^\prime(\varphi;\ParV,\ParV^\prime) &:=&
\frac1{|\Lambda_n|} \ism[\Lambda_n]{ \Big| e^{-\VIPar{x}{\varphi}{\ParV}} -e^{-\VIPar{x}{\varphi}{\ParV^\prime} }
\Big| } \\
W_{2,\Lambda_n}^\prime(\varphi;\ParV,\ParV^\prime) &:=&\frac1{|\Lambda_n|}  \sum_{x \in \varphi_{\Lambda_n}} \Big| \VIPar{x}{\varphi\setminus x}{\ParV}-
\VIPar{x}{\varphi \setminus x}{\ParV^\prime} \Big|. \\
\end{eqnarray*}
Let $k\geq 1$ and let $\ParV, \ParV^\prime \in \SpPar$ such that $||\ParV-\ParV^\prime||\leq \frac1k$, then under the assumption \textbf{[C1]} and from Theorem~\ref{thmNguyen} and Corollary~\ref{cor-CampGlotz}, we have $P_{\ParVT}-$almost surely as $n \to +\infty$
\begin{eqnarray*}
W_{1,\Lambda_n}^\prime(\Phi;\ParV,\ParV^\prime) &\longrightarrow& \Esp\left(
\left| e^{- V^{}\left( 0|\Phi; \ParV \right)} 
-e^{- V^{}\left( 0|\Phi; \ParV^\prime \right)} 
\right| \right) \\
W_{2,\Lambda_n}^\prime(\Phi;\ParV,\ParV^\prime) &\longrightarrow& \Esp\left(
\left|  V^{}\left( 0|\Phi; \ParV \right)-  V^{}\left( 0|\Phi; \ParV^\prime \right)    \right| e^{- V^{}\left( 0|\Phi; \ParVT \right)} 
\right) 
\end{eqnarray*}
Under Assumption \textbf{[C4]}, one may apply the mean value theorem in ${\RR[p]}$ as follows: there exist $\Vect{\xi}^{(1)},\ldots, \Vect{\xi}^{(p)} \in \prod_{j=1}^p \left[ \min(\theta_j,\theta_j^\prime),\max(\theta_j,\theta_j^\prime)\right]$ such that for all $\varphi\in \cSp$
$$
e^{- V^{}\left( 0|\varphi; \ParV \right) } -e^{ - V^{}\left( 0|\varphi; \ParV^\prime \right)}  = \sum_{j=1}^p \left(\theta_j-\theta_j^\prime \right)  {\frac{\partial V^{}}{\partial \theta_{j}}}  \left( 0|\varphi; {\Vect{\xi}}^{(j)} \right)
e^{ -  V^{}\left( 0|\varphi; \Vect{\xi}^{(j)} \right) }.
$$
This leads, under Assumption \textbf{[C4]}, to the following inequality
\begin{eqnarray*}
\Esp\left( \left|e^{- V^{}\left( 0|\Phi; \ParV \right) } -e^{ - V^{}\left( 0|\Phi; \ParV^\prime \right)} \right|\right)^2 & \leq& \Esp\left( \left|e^{- V^{}\left( 0|\Phi; \ParV \right) } -e^{ - V^{}\left( 0|\Phi; \ParV^\prime \right)} \right|^2\right) \\
&\leq& \Esp\left(||\ParV-\ParV^\prime||^2 \sum_{j=1}^p  \left| {\frac{\partial V^{}}{\partial \theta_{j}}}  \left( 0|\Phi; {\Vect{\xi}}^{(j)} \right)
e^{ -  V^{}\left( 0|\Phi; \Vect{\xi}^{(j)} \right) }\right|^2.
\right) \\
&\leq& \left(\frac1k\right)^2 \gamma_1^2, 
\end{eqnarray*}
with $\gamma_1:= \Esp\left(\sum_{j=1}^p
\max_{\ParV\in\SpPar}\left| {\frac{\partial V^{}}{\partial \theta_{j}}}  \left( 0|\Phi; \ParV \right)
e^{ -  V^{}\left( 0|\Phi; \ParV \right) }\right|^2
\right)<+\infty$. In such a way, one may also prove that
$$
\Esp\left(
\left|  V^{}\left( 0|\Phi; \ParV \right)-  V^{}\left( 0|\Phi; \ParV^\prime \right)    \right| e^{- V^{}\left( 0|\Phi; \ParVT \right)} 
\right)^2 \leq \left(\frac1k\right)^2 \gamma_2^2, 
$$
with $\gamma_2:= \Esp\left(\sum_{j=1}^p
\max_{\ParV\in\SpPar}\left| {\frac{\partial V^{}}{\partial \theta_{j}}}  \left( 0|\Phi; \ParV \right)
e^{ -  V^{}\left( 0|\Phi; \ParVT \right) }\right|^2
\right)$. Hence, for all $k\geq 1$ and for all $\ParV,\ParV^\prime\in \SpPar$ such that $||\ParV-\ParV^\prime||\leq\frac1k$ there exists $n_0(k) \geq 1$ such that for all $n\geq n_0(k)$, we have
$$
W_{1,\Lambda_n}^\prime \left( \varphi;\ParV,\ParV^\prime\right) \leq \frac2k \gamma_1 \qquad \mbox{ and }\qquad  W_{2,\Lambda_n}^\prime \left( \varphi;\ParV,\ParV^\prime\right) \leq \frac2k \gamma_2, \mbox{ for } P_{\ParVT}-\mbox{a.e. } \varphi.
$$
Since $\gamma_1$ and $\gamma_2$ are independent of $\ParV$ and $\ParV^\prime$, we have for all $n\geq n_0(k)$
$$
W_n\left(\varphi,\frac1k\right) \leq W_{1,n}\left(\varphi,\frac1k\right)+W_{2,n}\left(\varphi,\frac1k\right) \leq \frac2k\left( \gamma_1+\gamma_2\right):=\frac ck, \mbox{ for } P_{\ParVT}-\mbox{a.e. } \varphi.
$$ 
Finally, since
$$
\limsup_{n\to +\infty} \left\{ W_n\left(\varphi, \frac1k \right) \geq \frac ck \right\} = \bigcap_{m \in \NN}\bigcup_{n \geq m } \left\{ W_{n}\left( \varphi,\frac1k \right) \geq \frac ck  \right\} \subset
\bigcup_{n \geq n_0(k)} \left\{ W_{n}\left(\varphi, \frac1k \right) \geq \frac ck  \right\}
$$
for $P_{\ParVT}-$a.e. $\varphi$, the expected result (\ref{modCont}) is proved.

\noindent\textit{Conclusion step.}
The Steps 1, 2 and 3 ensure the fact that we can apply Property 3.6 of \cite{Guyon92} which asserts the almost sure convergence for minimum contrast estimators.


\subsection{Proof of Theorem~\ref{thm-norm}}

\noindent \textit{Step 1. Asymptotic normality of \ensuremath{\Vect{U}_n^{(1)}(\Phi;\ParVT)}}

The aim is to prove the following convergence in distribution as $n \to +\infty$
\begin{equation} \label{Un1Loi}
{}|\Lambda_n|^{1/2} \; \ensuremath{\Vect{U}_n^{(1)}(\Phi;\ParVT)} \rightarrow \mathcal{N} \left( 0,  \Mat{\Sigma}(\Dt,\ParVT)\right)
\end{equation}
where the matrix $ \Mat{\Sigma}(\Dt,\ParVT)$ is defined by (\ref{eq-defSig}).

The idea is to apply to $\ensuremath{\Vect{U}_n^{(1)}(\Phi;\ParVT)}$ a central limit theorem obtained by \cite{Jensen94}, Theorem~2.1. The following conditions have to be fulfilled to apply this result. For all $j=1,\ldots,p$
\begin{itemize}
\item[$(i)$] For all $i \in \ZZ[2]$, $\Esp \left(
\left(\ensuremath{\Vect{LPL}_{\Dom[i]}^{(1)} \left( \Phi ; \ParVT  \right)} \right)_j | \Phi_{\Dom[i]^c} \right)=0.$
\item[$(ii)$] For all $i\in \ZZ[2]$,  $\Esp \left( \left| \left( \ensuremath{\Vect{LPL}_{\Dom[i]}^{(1)} \left( \Phi ; \ParVT  \right)} \right)_j \right|^3 \right) <+\infty.$
\item[$(iii)$] The matrix $\Var \left( |\Lambda_n|^{1/2}  \ensuremath{\Vect{U}_n^{(1)}(\Phi;\ParVT)}\right)$ converges to the matrix $\Mat{\Sigma}(\Dt,\ParVT)$.
\end{itemize}

\noindent \underline{Condition $(i)$}~: From the stationarity of the process, it is sufficient to prove that 
$$\Esp \left( \left(\ensuremath{\Vect{LPL}_{\Delta_0}^{(1)} \left( \Phi ; \ParVT  \right)} \right)_j | \Phi_{\Delta_0^c} \right)=0.$$
Recall that for any configuration $\varphi$
\begin{equation} \label{goj}
 \left(\ensuremath{\Vect{LPL}_{\Delta_0}^{(1)} \left( \varphi ; \ParVT  \right)}\right)_j= -\ism[\Delta_0]{ \dVIPar{j}{x}{\varphi}{\ParVT}  e^{-\VIPar{x}{\varphi}{\ParVT}}} + 
\int_{\Delta_0}\dVIPar{j}{x}{\varphi \setminus x}{\ParVT} \varphi(dx).
\end{equation}
Denote respectively by $G_1(\varphi)$ and $G_2(\varphi)$ the first and the second right-hand term of~(\ref{goj}) and by $E_i=\Esp\left( G_i(\Phi) | \Phi_{\Delta_0^c}=\varphi_{\Delta_0^c} \right)$. Let us define for any $\varphi$, the measure $\mu_\varphi:=\sum_{x\in \varphi} \delta_x$.
From the definition of Gibbs point processes,
$$
E_2 = \frac1{Z_{\Delta_0}(\varphi_{\Delta_0^c})} \int_{\Omega_{\Delta_0}} \pi_{\Delta_0}(d\varphi_{\Delta_0}) 
\int_{\RR[2]} \mu_{\varphi_{\Delta_0}}(dx)\mathbf{1}_{\Delta_0}(x) \dVIPar{j}{x}{\varphi\setminus x}{\ParVT} e^{-
V_{\Delta_0}\left(\varphi;\ParVT\right)}.
$$
Since $\pi$ is a Poisson process, 
$$\int_{\cSp_{\Delta_0}} \pi_{\Delta_0}(d\varphi_{\Delta_0})f(\varphi)
=\int_{\cSp_{\Delta_0}} \pi_{\Delta_0}(d\varphi_{\Delta_0}) 
\int_{\cSp_{\Delta_0^c}} \pi_{\Delta_0^c}(d\varphi^\prime_{\Delta_0^c})f(\varphi)$$
and therefore, by introducing $\psi:=\varphi_{\Delta_0} \cup \varphi_{\Delta_0^c}^\prime$
$$
E_2=\frac1{Z_{\Delta_0}(\varphi_{\Delta_0^c})} \int_{\cSp} \pi(d\psi) 
\int_{\RR[2]} \mu_\psi(dx) \mathbf{1}_{\Delta_0}(x) 
\dVIPar{j}{x}{\psi_{\Delta_0} \cup \varphi_{\Delta_0^c} \setminus x}{\ParVT} 
e^{-V_{\Delta_0}\left(\psi_{\Delta_0} \cup \varphi_{\Delta_0^c};\ParVT\right)}.
$$
Now, from Campbell Theorem (applied to the Poisson measure $\pi$)
$$E_2=\frac1{Z_{\Delta_0}(\varphi_{\Delta_0^c})} \int_{\Delta_0}  dx \int_{\cSp} \pi_{x}^!(d\psi) 
\dVIPar{j}{x}{\psi_{\Delta_0} \cup \varphi_{\Delta_0^c}}{\ParVT}
e^{-V_{\Delta_0}\left(x\cup \psi_{\Delta_0} \cup \varphi_{\Delta_0^c};\ParVT\right)},
$$
where $\pi_x^!$ stands for the reduced Palm distribution of the Poisson point process. Since from Slivnyak-Mecke Theorem (see {\it e.g.} \cite{B-MolWaa03}), $\pi=\pi_x^!$, one can obtain
\begin{eqnarray}
E_2 &=& \frac1{Z_{\Delta_0}(\varphi_{\Delta_0^c})} \int_{\cSp} \pi(d\psi) 
\int_{\Delta_0}  dx \; 
\dVIPar{j}{x}{\psi_{\Delta_0} \cup \varphi_{\Delta_0^c}}{\ParVT}
e^{-V_{\Delta_0}\left(x\cup\psi_{\Delta_0} \cup \varphi_{\Delta_0^c};\ParVT\right)}
 \nonumber \\
&=& \frac1{Z_{\Delta_0}(\varphi_{\Delta_0^c})} \int_{\cSp_{\Delta_0}} \pi_{\Delta_0}(d\varphi_{\Delta_0})
\int_{\Delta_0} dx \dVIPar{j}{x}{\varphi}{\ParVT}
e^{- \VIPar{x}{\varphi}{\ParVT}} 
e^{-V_{\Delta_0}\left(\varphi;\ParVT\right)}
\nonumber \\
&=& -E_1 \nonumber 
\end{eqnarray}

\noindent \underline{Condition $(ii)$}~: 
For any bounded domain $\Delta$ one may write for~$j=1,\ldots,p$
$$
\left| \left(\ensuremath{\Vect{LPL}_{\Delta}^{(1)} \left( \Phi ; \ParVT  \right)} \right)_j\right|^3 \leq 4 \left| \ism[\Delta] {\dVIPar{j}{x}{\Phi}{\ParVT} e^{- \VIPar{x}{\Phi}{\ParVT} }}
\right|^3 +
4 \left| \sum_{x \in \varphi_{\Delta}} \dVIPar{j}{x}{\Phi \setminus x}{\ParVT}
\right|^3.
$$
The assumption \textbf{[N1]} ensures  the integrability of the first right-hand term. For the second one, note that
\begin{eqnarray*}
T_2&:=&\left| \sum_{x \in \Phi_{\Delta}} \dVIPar{j}{x}{\varphi \setminus x}{\ParVT}
\right|^3 \\
&\leq& \mathop{\sum_{x_1, x_2, x_3 \in \varphi_\Delta}}_{x_1\neq x_1, x_2\neq x_3 ,x_2\neq x_3}
\left|\dVIPar{j}{x_1}{\varphi \setminus x_1}{\ParVT}\right| \left|\dVIPar{j}{x_2}{\varphi \setminus x_2}{\ParVT}\right| \left|\dVIPar{j}{x_3}{\varphi \setminus x_3}{\ParVT}\right| \\
&&+ 3 \sum_{x_1, x_2 \in \varphi_\Delta,x_1\neq x_2}
\left|\dVIPar{j}{x_1}{\varphi \setminus x_1}{\ParVT}\right|^2
\left|\dVIPar{j}{x_2}{\varphi \setminus x_2}{\ParVT}\right| \\
&&+\sum_{x_1 \in \varphi_\Delta}
\left|\dVIPar{j}{x_2}{\varphi \setminus x_1}{\ParVT}\right|^3.
\end{eqnarray*}
The result is obtained by using the assumption \textbf{[N1]} and  iterated versions of Corollary~\ref{cor-CampGlotz}.\\

\noindent \underline{Condition $(iii)$:} let us start by noting that from the assumption \textbf{[Mod-L]}, the vector $\ensuremath{\Vect{LPL}_{\Dom[i]}^{(1)} \left( \varphi ; \ParVT  \right)}$ depends only on $\varphi_{\Dom[j]}$ for $j\in \BDt[i]{1}$. Let $E_{i,j}:= \Esp \left( \ensuremath{\Vect{LPL}_{\Dom[i]}^{(1)} \left( \Phi ; \ParVT  \right)} \tr{\ensuremath{\Vect{LPL}_{\Dom[j]}^{(1)} \left( \Phi ; \ParVT  \right)}} \right)$. Based on our definitions, we have
\begin{eqnarray*}
\Var \left( |\Lambda_n|^{1/2}  \ensuremath{\Vect{U}_n^{(1)}(\Phi;\ParVT)} \right)
&=& |\Lambda_n|^{-1} \Var \left( \sum_{i \in I_n}\ensuremath{\Vect{LPL}_{\Dom[i]}^{(1)} \left( \Phi ; \ParVT  \right)} \right) \\
&=&|\Lambda_n|^{-1} \sum_{i,j \in I_n} E_{i,j}  \\
&=& |\Lambda_n|^{-1} \sum_{i \in I_n} \left( 
\sum_{j\in I_n \cap \BDt[i]{1}}   E_{i,j}  + \sum_{j\in I_n \cap \BDt[i]{1}^c}  E_{i,j} \right).
\end{eqnarray*}
Let $j\in I_n \cap \BDt[i]{1}^c$, since \ensuremath{\Vect{LPL}_{\Delta_i}^{(1)} \left( \varphi ; \ParVT  \right)} is a measurable function of $\varphi_{\Delta_i^c}$, we have by using condition~(i):
\begin{eqnarray}
\Esp \left(  \ensuremath{\Vect{LPL}_{\Dom[i]}^{(1)} \left( \Phi ; \ParVT  \right)} \tr{\ensuremath{\Vect{LPL}_{\Dom[j]}^{(1)} \left( \Phi ; \ParVT  \right)}} \right) &=& \Esp \left( \Esp \left( \ensuremath{\Vect{LPL}_{\Dom[i]}^{(1)} \left( \Phi ; \ParVT  \right)} \tr{\ensuremath{\Vect{LPL}_{\Dom[j]}^{(1)} \left( \Phi ; \ParVT  \right)}} | \Phi_{\Dom[i]^c}  \right)\right) \nonumber \\
&=& \Esp\left(\Esp \left( \ensuremath{\Vect{LPL}_{\Dom[i]}^{(1)} \left( \Phi ; \ParVT  \right)} | \Phi_{\Dom[i]^c}
 \right)  \tr{\ensuremath{\Vect{LPL}_{\Dom[j]}^{(1)} \left( \Phi ; \ParVT  \right)}}\right) \nonumber \\
&=& 0 \nonumber
\end{eqnarray}
Denote by $\widetilde{I}_n$ the following set
$$
\widetilde{I}_n = I_n \cap \left( \cup_{i\in \partial I_n} \BDt[i]{1} \right).
$$
We now obtain
\begin{eqnarray*}
\Var \left( |\Lambda_n|^{1/2}  \ensuremath{\Vect{U}_n^{(1)}(\Phi;\ParVT)} \right)
&=&  |\Lambda_n|^{-1} \sum_{i \in I_n} \sum_{j\in I_n \cap \BDt[i]{1}}   E_{i,j} \\
&=& |\Lambda_n|^{-1} \left( 
\sum_{i \in I_n\setminus\widetilde{I}_n} \sum_{j\in I_n \cap \BDt[i]{1}} E_{i,j} + \sum_{i \in \widetilde{I}_n} \sum_{j\in I_n \cap \BDt[i]{1}} E_{i,j}
\right)
\end{eqnarray*}
Using the stationarity and the definition of the domain $\Lambda_n$, one obtains
$$
{}|\Lambda_n|^{-1} 
\sum_{i \in I_n\setminus\widetilde{I}_n} \sum_{j\in I_n \cap \BDt[i]{1}} E_{i,j} = |\Lambda_n|^{-1}  |I_n\setminus\widetilde{I}_n| \sum_{j\in \BDt[0]{1}} E_{0,j} \to \Mat{\Sigma}(\Dt,\ParVT) \quad \mbox{ as } n\to +\infty
$$
and
$$
{}|\Lambda_n|^{-1} \left| \sum_{i \in \widetilde{I}_n} \sum_{j\in I_n \cap \BDt[i]{\cD[\frac{D}{\Dt}]}} E_{i,j}
\right| \leq |\Lambda_n|^{-1} |\widetilde{I}_n| \sum_{j\in \BDt[0]{1}} |E_{0,j}| \to 0  \quad \mbox{ as } n\to +\infty.
$$
Hence as $n\to +\infty$  
\begin{eqnarray}
\Var \left( |\Lambda_n|^{1/2}  \ensuremath{\Vect{U}_n^{(1)}(\Phi;\ParVT)} \right) &=&  |\Lambda_n|^{-1} \sum_{i \in I_n}  \sum_{j\in I_n \cap \BDt[i]{1}}   E_{i,j} \nonumber \\ 
&\stackrel{n\to +\infty}{\longrightarrow}& \underbrace{|I_n| |\Lambda_n|^{-1}}_{\Dt^{-2}}\sum_{k\in \BDt[0]{1}} E_{0,k} = \Mat{\Sigma}(\Dt,\ParVT).\label{eq-convVar} 
\end{eqnarray}


\noindent \textit{Step 2. Domination of $ \ensuremath{\Mat{U}_n^{(2)}(\Phi;\ParV)}$ in a neighborhood of $\ParVT$ and convergence of $\ensuremath{\Mat{U}_n^{(2)}(\Phi;\ParVT)}$}
Let $j,k=1,\ldots,p$, recall that $\left( \ensuremath{\Mat{U}_n^{(2)}(\varphi;\ParV)}\right)_{j,k}$ is defined in a neighborhood $\Vois(\ParVT)$ of $\ParVT$ for any configuration $\varphi$ by
\begin{eqnarray}
\left(\ensuremath{\Mat{U}_n^{(2)}(\varphi;\ParV)}
\right)_{j,k} &=& -\frac1{|\Lambda_n|} \ism[\Lambda_n ]{
\ddVIPar{j}{k}{x}{\varphi}{\ParV} 
\exp\left( - \VIPar{x}{\varphi}{\ParV}\right) } \nonumber \\ 
&& +  
\frac1{|\Lambda_n|} \ism[\Lambda_n]{
\dVIPar{j}{x}{\varphi}{\ParV} \dVIPar{k}{x}{\varphi}{\ParV}
\exp\left( - \VIPar{x}{\varphi}{\ParV}\right)} \nonumber \\
&& + \frac1{|\Lambda_n|} \sum_{x \in \varphi_{\Lambda_n}} 
\ddVIPar{j}{k}{x}{\varphi\setminus x}{\ParV} .
\end{eqnarray}
Under the assumption \textbf{[N1]} and \textbf{[N2]}, from Theorem~\ref{thmNguyen} (\cite{Nguyen79}) and from Corollary~\ref{cor-CampGlotz}, there exists $n_0 \in \NN$ such that for all $n\geq n_0$
\begin{eqnarray*}
\left|\left(\ensuremath{\Mat{U}_n^{(2)}(\varphi;\ParV)} \right)_{j,k}\right| &\leq& 2 \Esp\left( \left(\left| \ddVIPar{j}{k}{0}{\Phi}{\ParV} 
\right| 
+\left| \dVIPar{j}{0}{\Phi}{\ParV}  \dVIPar{k}{0}{\Phi}{\ParV}  
\right|\right)
e^{-\VIPar{0}{\Phi}{\ParV}}
\right)\\
&& +2\times \Esp\left( \left| \ddVIPar{j}{k}{0}{\Phi}{\ParV} 
\right| e^{-\VIPar{0}{\Phi}{\ParVT}}
\right)
\end{eqnarray*}
Note that from~Theorem~\ref{thmNguyen} (\cite{Nguyen79}), $\ensuremath{\Mat{U}_n^{(2)}(\cdot;\ParVT)}$ converges almost surely as $n \to +\infty$ towards $\Mat{U}^{(2)}(\ParVT)$ defined by~(\ref{def-U2}). Note that $\Mat{U}^{(2)}(\ParVT)$ is a symmetric positive matrix since for all $\Vect{y}\in \RR[p]$
$$
\tr{\Vect{y}} \Mat{U}^{(2)}(\ParVT) \Vect{y} = 
\Esp\left(
\left(\tr{\Vect{y}} \Vect{V}^{(1)}(0|\Phi;\ParVT)\right)^2 e^{-\VIPar{0}{\Phi}{\ParVT}}
\right) \geq 0,
$$
where for $j=1,\ldots,p$, $\varphi\in \cSp$ and for $\ParV\in \Vois(\ParVT)$ $\left(\Vect{V}^{(1)}(x|\varphi;\ParVT)\right)_j :=\dVIPar{j}{x}{\varphi}{\ParV}$ and it is a definite matrix under the assumption~\textbf{[N3]}.

%
\noindent \textit{Conclusion Step} 
Under the assumptions \textbf{[Mod]} and \textbf{[Ident]}, and using Steps 1 and 2, one can apply a classical result concerning asymptotic normality for minimum contrast estimators 
{\it e.g.} Proposition 3.7 of \cite{Guyon92} in order to obtain~(\ref{convLoiMPLE1}). \\

It remains to prove (\ref{convLoiMPLE2}). This may de done in two different steps. The first one consists in verifying the positive definiteness of the matrix $\Mat{\Sigma}(\Dt,\ParVT)$. The proof is strictly similar to the one of \cite{A-BilCoeDro08} (p.~261) except that the assumption [\textbf{SDP}] is now simply replaced by the more general one assumption \textbf{[N4]}. Now, the convergence in probability of $\ensuremath{\Estn{\Mat{\Sigma}}(\Phi;\Estn{\ParV}(\Phi))}$ towards $\Mat{\Sigma}(\Dt,\ParVT)$ is obtained by applying Proposition~9 of \cite{A-CoeLav10}.

\section*{Acknowledgements}

We are grateful to the referee and the associate editor for their comments which helped us in improving a previous version. 

\bibliographystyle{plainnat.bst}
\bibliography{mpleEJS}

\end{document}

%% file: cqlsInclude.tex
\usepackage{mathtools}
\usepackage{amssymb}

\newcommand{\cqlsbm}[1]{\mbox{\boldmath$#1$}}

\newcommand{\Var}{
        \mathbb{V}\mbox{ar}
}

\newcommand{\Mat}[1]{
  \underline{\cqlsbm{#1}}
}

\newcommand{\tr}[1]{
 {#1}^{\!  T}
}

\newcommand{\Vect}[1]{
  \cqlsbm{#1}
}

\newcommand{\Esp}{
  \cqlsbm{E}
}

\def\RR{\hbox{I\kern-.1em\hbox{R}}}
\def\NN{\hbox{I\kern-.1em\hbox{N}}}
\def\ZZ{\hbox{Z\kern-.3em\hbox{Z}}}
\def\PP{\hbox{P\kern-.8em\hbox{P}}}

%% file: testInclude.tex
\newlength{\hauteur}
\newlength{\profond}

%% file: mpleGenHAL2.bbl
\begin{thebibliography}{29}
\providecommand{\natexlab}[1]{#1}
\providecommand{\url}[1]{\texttt{#1}}
\expandafter\ifx\csname urlstyle\endcsname\relax
  \providecommand{\doi}[1]{doi: #1}\else
  \providecommand{\doi}{doi: \begingroup \urlstyle{rm}\Url}\fi

\bibitem[Baddeley and Turner(2000)]{A-BadTur00}
A.~Baddeley and R.~Turner.
\newblock {{P}ractical maximum pseudolikelihood for spatial point patterns
  (with discussion)}.
\newblock \emph{Australian and New Zealand Journal of Statistics}, 42:\penalty0
  283--322, 2000.

\bibitem[Bertin et~al.(1999{\natexlab{a}})Bertin, Billiot, and
  Drouilhet]{A-BerBilDro99}
E.~Bertin, J.-M. Billiot, and R.~Drouilhet.
\newblock Existence of ``{N}earest-{N}eighbour'' {G}ibbs {P}oint {M}odels.
\newblock \emph{Ann. Appl. Probab.}, 31:\penalty0 895--909, 1999{\natexlab{a}}.

\bibitem[Bertin et~al.(1999{\natexlab{b}})Bertin, Billiot, and Drouilhet]{BBD1}
E.~Bertin, J.-M. Billiot, and R.~Drouilhet.
\newblock Spatial {D}elaunay {G}ibbs {Point} {P}rocesses.
\newblock \emph{Stochastic {M}odels}, 15\penalty0 (2):\penalty0 181--199,
  1999{\natexlab{b}}.

\bibitem[Bertin et~al.(1999{\natexlab{c}})Bertin, Billiot, and Drouilhet]{BBD2}
E.~Bertin, J.-M. Billiot, and R.~Drouilhet.
\newblock $k$-{N}earest-{N}eighbour {G}ibbs {P}oint {P}rocesses.
\newblock \emph{Markov {P}rocesses and {R}elated {F}ields}, 5\penalty0
  (2):\penalty0 219--234, 1999{\natexlab{c}}.

\bibitem[Besag(1974)]{Besag74}
J.~Besag.
\newblock Spatial interaction and the statistical analysis of lattice system.
\newblock \emph{J. R. Statist. Soc. Ser. B}, 26:\penalty0 192--236, 1974.

\bibitem[Besag et~al.(1982)Besag, Milne, and Zachary]{Besag82}
J.~Besag, R.~Milne, and S.~Zachary.
\newblock Point process limits of lattice processes.
\newblock \emph{Ann. Appl. Prob.}, 19:\penalty0 210--216, 1982.

\bibitem[Billiot et~al.(2008)Billiot, Coeurjolly, and Drouilhet]{A-BilCoeDro08}
J.-M. Billiot, J.-F. Coeurjolly, and R.~Drouilhet.
\newblock {M}aximum pseudolikelihood estimator for exponential family models of
  marked {G}ibbs point processes.
\newblock \emph{Electronic Journal of Statistics}, 2:\penalty0 234--264, 2008.

\bibitem[Coeurjolly and Lavancier(2010)]{A-CoeLav10}
J.-F. Coeurjolly and F.~Lavancier.
\newblock Residuals for stationary marked {G}ibbs point processes.
\newblock \emph{submitted for publication}, 2010.
\newblock \texttt{http://hal.archives-ouvertes.fr/hal-00453102/fr/}.

\bibitem[Daley and Vere-Jones(1988)]{B-DalVer88}
D.~Daley and D.~Vere-Jones.
\newblock \emph{An introduction to the {T}heory of {P}oint {P}rocesses}.
\newblock Springer Verlag, New York, 1988.

\bibitem[Dereudre(2005)]{A-Der05}
D.~Dereudre.
\newblock {G}ibbs {D}elaunay tessellations with geometric hardcore condition.
\newblock \emph{J. Stat. Phys.}, 121\penalty0 (3-4):\penalty0 511--515, 2005.

\bibitem[Dereudre(2009)]{A-Der09}
D.~Dereudre.
\newblock {The existence of quermass-interaction processes for nonlocally
  stable interaction and nonbounded convex grains}.
\newblock \emph{Adv. in Appl. Probab}, 41\penalty0 (3):\penalty0 664--681,
  2009.

\bibitem[Dereudre and Lavancier(2009)]{A-DerLav09}
D.~Dereudre and F.~Lavancier.
\newblock {C}ampbell equilibrium equation and pseudo-likelihood estimation for
  non-hereditary {G}ibbs point processes.
\newblock \emph{Bernoulli}, 15\penalty0 (4):\penalty0 1368--1396, 2009.

\bibitem[Dereudre et~al.(2010)Dereudre, Drouilhet, and Georgii]{A-DerDroGeo09}
D.~Dereudre, R.~Drouilhet, and H.O. Georgii.
\newblock Existence of {G}ibbsian point processes with geometry-dependent
  interactions.
\newblock \emph{submitted}, 2010.
\newblock http://arxiv.org/abs/1003.2875.

\bibitem[Goulard et~al.(1996)Goulard, S\"arkk{\"a}, and
  Grabarnik]{A-GouSarGra96}
M.~Goulard, A.~S\"arkk{\"a}, and P.~Grabarnik.
\newblock {Parameter estimation for marked Gibbs point processes through the
  maximum pseudo-likelihood method}.
\newblock \emph{Scandinavian Journal of Statistics}, 23\penalty0 (3):\penalty0
  365--379, 1996.

\bibitem[Guyon(1992)]{Guyon92}
X.~Guyon.
\newblock \emph{Champs al\'eatoires sur un r\'eseau}.
\newblock Masson, Paris, 1992.

\bibitem[Illian et~al.(2008)Illian, H., and Stoyan]{B-IllPenSto08}
J.~Illian, A.~Penttinen H., and Stoyan.
\newblock \emph{{Statistical analysis and modelling of spatial point
  patterns}}.
\newblock Wiley-Interscience, 2008.

\bibitem[Jensen and Künsch(1994)]{Jensen94}
J.L. Jensen and H.R. Künsch.
\newblock On asymptotic normality of pseudo likelihood estimates of pairwise
  interaction processes.
\newblock \emph{Ann. Inst. Statist. Math.}, 46:\penalty0 475--486, 1994.

\bibitem[Jensen and M\o{}ller(1991)]{Jensen91}
J.L. Jensen and J.~M\o{}ller.
\newblock Pseudolikelihood for exponential family models of spatial point
  processes.
\newblock \emph{Ann. Appl. Probab.}, 1:\penalty0 445--461, 1991.

\bibitem[Kendall et~al.(1999)Kendall, Lieshout, and Baddeley]{A-KenLieBad99}
W.~S. Kendall, M.~N. M.~Van Lieshout, and A.~J. Baddeley.
\newblock Quermass-interaction processes conditions for stability.
\newblock \emph{Advances on applied probability}, 31:\penalty0 315--342, 1999.

\bibitem[Mase(1995)]{Mase95}
S.~Mase.
\newblock Consistency of maximum pseudo-likelihood estimator of continuous
  state space gibbsian process.
\newblock \emph{Ann. Appl. Probab.}, 5:\penalty0 603--612, 1995.

\bibitem[Mase(2000)]{Mase00}
S.~Mase.
\newblock Marked gibbs processes and asymptotic normality of maximum
  pseudo-likelihhod estimators.
\newblock \emph{Math. Nachr.}, 209:\penalty0 151--169, 2000.

\bibitem[M{\o}ller(2008)]{R-Mol08}
J.~M{\o}ller.
\newblock {P}arametric methods for spatial point processes.
\newblock Technical Report Research Report R-2008-04, Department of
  Mathematical Sciences, Aalborg University, 2008.

\bibitem[M\o{}ller and Waagepetersen(2003)]{B-MolWaa03}
J.~M\o{}ller and R.~Waagepetersen.
\newblock \emph{Statistical Inference and Simulation for Spatial Point
  Processes}.
\newblock Chapman and Hall/CRC, Boca Raton, 2003.

\bibitem[Nguyen and Zessin(1979)]{Nguyen79}
X.X. Nguyen and H.~Zessin.
\newblock Ergodic theorems for {S}patial {P}rocess.
\newblock \emph{Z. Wahrscheinlichkeitstheorie verw. Gebiete}, 48:\penalty0
  133--158, 1979.

\bibitem[Ogata and Tanemura(1981)]{A-OgaTan81}
Y.~Ogata and M.~Tanemura.
\newblock {Estimation of interaction potentials of spatial point patterns
  through the maximum likelihood procedure}.
\newblock \emph{Annals of the Institute of Statistical Mathematics},
  33\penalty0 (1):\penalty0 315--338, 1981.

\bibitem[Preston(1976)]{B-Pre76}
C.J. Preston.
\newblock \emph{Random fields}.
\newblock Springer Verlag, 1976.

\bibitem[Ruelle(1969)]{B-Rue69}
D.~Ruelle.
\newblock \emph{{S}tatistical {M}echanics}.
\newblock Benjamin, New York-Amsterdam, 1969.

\bibitem[Ruelle(1970)]{Ruelle70}
D.~Ruelle.
\newblock Superstable interactions in classical statistical mechanics.
\newblock \emph{Commun. Math. Phys.}, 18:\penalty0 127--159, 1970.

\bibitem[Stoyan et~al.(1987)Stoyan, Kendall, Mecke, and
  Ruschendorf]{B-StoKenMecRus87}
D.~Stoyan, W.S. Kendall, J.~Mecke, and L.~Ruschendorf.
\newblock \emph{{{S}tochastic geometry and its applications}}.
\newblock John Wiley and Sons, Chichester, 1987.

\end{thebibliography}
